\documentclass[11pt]{amsart}

\usepackage{amsmath}
\usepackage{amscd}
\usepackage{amssymb}
\usepackage{latexsym}
\usepackage{stmaryrd} 
\usepackage{mathbbol} 

\usepackage[arrow,curve,matrix,tips,2cell]{xy}
  \SelectTips{eu}{10} \UseTips
  \UseAllTwocells
\newtheorem{theorem}{Theorem}[section]
\newtheorem*{theorem*}{Theorem}
\newtheorem{lemma}[theorem]{Lemma}
\newtheorem*{lemma*}{Lemma}
\newtheorem{corollary}[theorem]{Corollary}
\newtheorem{proposition}[theorem]{Proposition}

\newtheorem{remark}[theorem]{Remark}
\newtheorem{definition}[theorem]{Definition}

\newcommand{\bgl}{\begin{equation}} 
\newcommand{\egl}{\end{equation}}
\newcommand{\bgloz}{\begin{equation*}} 
\newcommand{\egloz}{\end{equation*}}
\newcommand{\bgln}{\begin{eqnarray}} 
\newcommand{\egln}{\end{eqnarray}}
\newcommand{\bglnoz}{\begin{eqnarray*}} 
\newcommand{\eglnoz}{\end{eqnarray*}}
\newcommand{\btheo}{\begin{theorem}}
\newcommand{\etheo}{\end{theorem}}
\newcommand{\btheooz}{\begin{theorem*}}
\newcommand{\etheooz}{\end{theorem*}}
\newcommand{\blemma}{\begin{lemma}}
\newcommand{\elemma}{\end{lemma}}
\newcommand{\blemmaoz}{\begin{lemma*}}
\newcommand{\elemmaoz}{\end{lemma*}}
\newcommand{\bproof}{\begin{proof}}
\newcommand{\eproof}{\end{proof}}
\newcommand{\bbew}{\begin{beweis}}
\newcommand{\ebew}{\end{beweis}}
\newcommand{\bremark}{\begin{remark}\em}
\newcommand{\eremark}{\end{remark}}
\newcommand{\bdefin}{\begin{definition}}
\newcommand{\edefin}{\end{definition}}
\newcommand{\bprop}{\begin{proposition}}
\newcommand{\eprop}{\end{proposition}}
\newcommand{\bcor}{\begin{corollary}}
\newcommand{\ecor}{\end{corollary}}
\newcommand{\bfa}{\begin{cases}} 
\newcommand{\efa}{\end{cases}}
\newcommand{\cE}{\mathcal E}
\newcommand{\cF}{\mathcal F}
\newcommand{\cG}{\mathcal G}
\newcommand{\cH}{\mathcal H}

\newcommand{\cJ}{\mathcal J}
\newcommand{\cK}{\mathcal K}
\newcommand{\cL}{\mathcal L}

\newcommand{\cO}{\mathcal O}

\newcommand{\cU}{\mathcal U}

\def\Cz{\mathbb{C}}
\def\Fz{\mathbb{F}}

\def\Nz{\mathbb{N}}

\def\Zz{\mathbb{Z}}

\def\1z{\mathbb{1}}
\newcommand{\fA}{\mathfrak A}

\newcommand{\fJ}{\mathfrak J}

\newcommand{\fP}{\mathfrak P}

\newcommand{\an}[1]{``#1''} 
\newcommand{\ti}{\tilde}

\newcommand{\lori}{\longrightarrow}

\newcommand{\ma}{\mapsto} 
\newcommand{\mafr}{\mapsfrom} 
\newcommand\into{\hookrightarrow} 
\newcommand{\Rarr}{\Rightarrow} 
\newcommand{\Larr}{\Leftarrow} 
\newcommand{\LRarr}{\Leftrightarrow} 

\newcommand{\ve}{\varepsilon}

\def\SEMI{\mbox{$\times\kern-2pt\vrule height5pt width.6pt \kern3pt $}}

\newcommand{\PAut}{{\rm PAut}\,}

\newcommand{\Spec}{{\rm Spec\,}} 

\newcommand{\id}{{\rm id}}
\newcommand{\alg}{{\rm alg}}

\newcommand{\Ind}{\mathrm{ Ind}\,}

\newcommand{\Ad}{{\rm Ad\,}}

\renewcommand{\ker}{{\rm ker}\,}

\newcommand{\reg}{^\times} 
\newcommand{\tr}{\operatorname{tr}} 
\newcommand{\lspan}{{\rm span}} 
\newcommand{\clspan}{\overline{\lspan}} 
\newcommand{\abs}[1]{\lvert#1\rvert} 
\newcommand{\norm}[1]{\left\|#1\right\|} 
\newcommand{\defeq}{\mathrel{:=}} 
\newcommand{\eqdef}{\mathrel{=:}} 

\newcommand{\dop}{\text{: }} 
\newcommand{\falls}{\text{ if }} 
\newcommand{\sonst}{\text{ else}} 
\newcommand{\fa}{\text{ for all }} 
\newcommand{\e}[1]{e_{\left[#1\right]}} 
\newcommand{\E}[1]{E_{\left[#1\right]}} 
\newcommand{\rta}{\rtimes^a} 
\newcommand{\dsupp}{\text{d-supp}}

\newcommand{\Proj}{{\rm Proj}\,}

\newcommand{\lge}{\left\{} 
\newcommand{\rge}{\right\}} 
\newcommand{\lru}{\left(} 
\newcommand{\rru}{\right)} 
\newcommand{\leck}{\left[} 
\newcommand{\reck}{\right]} 
\newcommand{\lsp}{\left\langle} 
\newcommand{\rsp}{\right\rangle} 
\newcommand{\rukl}[1]{\lru #1 \rru} 
\newcommand{\eckl}[1]{\leck #1 \reck} 
\newcommand{\gekl}[1]{\lge #1 \rge} 
\newcommand{\spkl}[1]{\lsp #1 \rsp} 
\newcommand{\menge}[2]{\gekl{ #1 \dop #2 }} 
\newcommand{\ping}{P \subseteq G}
\newcommand{\aalphaP}{a_{(\alpha \vert_P)}}
\newcommand{\rtaas}{\rta_{\alpha,s}} 
\newcommand{\rtaar}{\rta_{\alpha,r}} 
\begin{document}

\title[Nuclearity of semigroup C*-algebras]{Nuclearity of semigroup C*-algebras and the connection to amenability}

\author{Xin Li}

\address{Xin Li, Department of Mathematics, Westf{\"a}lische Wilhelms-Universit{\"a}t M{\"u}nster, Einsteinstra{\ss}e 62, 48149 M{\"u}nster, Germany}
\email{xinli.math@uni-muenster.de}



\thanks{\scriptsize{Research supported by the ERC through AdG 267079.}}

\begin{abstract}
We study C*-algebras associated with subsemigroups of groups. For a large class of such semigroups including positive cones in quasi-lattice ordered groups and left Ore semigroups, we describe the corresponding semigroup C*-algebras as C*-algebras of inverse semigroups, groupoid C*-algebras and full corners in associated group crossed products. These descriptions allow us to characterize nuclearity of semigroup C*-algebras in terms of faithfulness of left regular representations and amenability of group actions. Moreover, we also determine when boundary quotients of semigroup C*-algebras are UCT Kirchberg algebras. This leads to a unified approach to Cuntz algebras and ring C*-algebras. 
\end{abstract}

\maketitle


\setlength{\parindent}{0pt} \setlength{\parskip}{0.5cm}

\section{Introduction}

We continue the project started in \cite{Li2} about C*-algebras associated with semigroups. The study of such semigroup C*-algebras goes back to L. Coburn (\cite{Co1}, \cite{Co2}) and was continued in for example \cite{Dou}, \cite{Mur1}, \cite{Mur2}, \cite{Mur3} and \cite{Mur4}. While there is a canonical reduced version, namely the C*-algebra generated by the left regular representation of the (left cancellative) semigroup, G. Murphy showed in \cite{Mur4} that the most obvious candidate for the full semigroup C*-algebra is intractable even for very simple (for instance abelian) semigroups. So one of the main difficulties was to find a good full version of semigroup C*-algebras, given by generators and relations, which could be viewed as the analogue of full group C*-algebras. 

One big step forward was \cite{Ni1}. A. Nica's idea was to define full semigroup C*-algebras using not only the obvious relations as in \cite{Mur3} but also additional ones reflecting the (right) ideal structure of the semigroup. This modification leads to interesting C*-algebras which can be analyzed and which exhibit good properties. However, A. Nica did not explicitly mention ideals of semigroups. Instead, he restricted his analysis to positive cones in quasi-lattice ordered groups which have a very simple ideal structure.

A. Nica's ideas have been taken up by M. Laca in collaboration with I. Raeburn and J. Crisp (\cite{La-Rae}, \cite{La1}, \cite{Cr-La1}, \cite{Cr-La2}). They studied the question when the left regular representation from the full to the reduced semigroup C*-algebra is faithful, and they described induced ideals of semigroup C*-algebras. However, the question when semigroup C*-algebras are nuclear was left untouched, and the connection between nuclearity and faithfulness of the left regular representation remained mysterious.

Recently, new examples of C*-algebras arising from number theory (\cite{Cun}, \cite{Cu-Li1}, \cite{Li1}, \cite{Cu-Li2}) have motivated the author to generalize A. Nica's work. For semigroups associated with number theoretic rings, the restriction to positive cones of quasi-lattice ordered groups corresponds to only considering principal ideal domains -- a restriction which, especially for rings from algebraic number theory, would exclude all the interesting examples. Making explicit use of the ideal structure of semigroups, the author was able to extend A. Nica's construction to arbitrary left cancellative semigroups in \cite{Li2}. The same construction was introduced independently in \cite{C-D-L} for particular examples of number theoretic interest. In general, it turns out that the full semigroup C*-algebras still have good properties. For instance, it is shown in \cite{Li2} and also \cite{Nor} that they are well-suited for studying amenability of semigroups. However, amenability is a strong assumption which interesting examples fail to have. One of the most striking examples is probably the $n$-fold free product $\Nz_0^{*n}$ of the natural numbers. This example is due to A. Nica, and he observed that it is closely related to the Cuntz algebra $\cO_n$.

A closely related topic is the theory of semigroup crossed products (by endomorphisms). One of the most important ideas in the analysis of semigroup crossed products is the idea of dilation. It already goes back to J. Cuntz in his work on the Cuntz algebras. This dilation theory has then been fully developed, in the probably most general setting, by M. Laca in \cite{La2}. He shows that one can use inductive limit procedures to dilate isometries to unitaries and endomorphisms to automorphisms so that in the end, semigroup crossed products can be embedded as full corners into group crossed products. This means that questions about semigroup crossed products translate into questions about group crossed products which have already been intensively studied. However, this dilation theory as described here only works for left Ore semigroups, and the question remains what to do for semigroups like the free product $\Nz_0^{*n}$.

Now, in the present paper, our main observation is that for semigroup C*-algebras in the sense of \cite{Ni1} or \cite{Li2}, the left Ore condition is not essential for embedding semigroup C*-algebras as full corners into group crossed products.

More precisely, for a subsemigroup $P$ of a group $G$, we show that under two conditions, the full and reduced semigroup C*-algebras of $P$ embed as full corners into full and reduced crossed products by $G$. The underlying (C*-)dynamical system is the same for both the full and reduced version. It is in a canonical way built out of the inclusion $\ping$ and a distinguished commutative subalgebra of the semigroup C*-algebras. The two conditions we have to impose are that the constructible right ideals of $P$ are independent and that $\ping$ satisfies the so-called Toeplitz condition. The first condition was introduced in \cite{Li2} and guarantees that the canonical commutative subalgebras of the full and reduced semigroup C*-algebras coincide. This condition also plays a crucial role in \cite{C-E-L1}. The second condition is new. It basically says that the procedure of compressing operators on $\ell^2(G)$ to $\ell^2(P)$ is well-behaved. We show that this condition is satisfied in typical examples. In particular, it holds for positive cones in quasi-lattice ordered groups and left Ore semigroups. Our main point is that we do not need the left Ore condition, only the two conditions described above. In order to embed full semigroup C*-algebras as full corners into group crossed products, the idea is to write both semigroup C*-algebras and the group crossed products into which we would like to embed as groupoid C*-algebras. The underlying groupoids are equivalent more or less by construction, so that we can use the observation by \cite{M-R-W} that equivalence of groupoids give rise to explicit imprimitivity bimodules of the corresponding groupoid C*-algebras. This result allows us to show that certain universal norms coincide. We point out that we work with the full version of semigroup C*-algebras introduced in \S~3 in \cite{Li2}.

As an application, we give equivalent characterizations for nuclearity of semigroup C*-algebras. For instance, we see that nuclearity can be expressed in terms of amenability of group actions. Moreover, nuclearity of semigroup C*-algebras implies faithfulness of the corresponding left regular representations.

In addition, we extend existing results about induced ideals and boundary actions from the quasi-lattice ordered case to our more general setting. This leads to a unified approach to specific constructions like Cuntz algebras or ring C*-algebras. As a second application of our main observation, we obtain a general explanation why these examples are UCT Kirchberg algebras.

A third application of our main observation is presented in \cite{C-E-L2} which constitutes a vast generalization of the K-theoretic results in \cite{C-E-L1}.

The present paper is structured as follows:

In a first preliminary section, we describe the setting (\S~\ref{setting}) and analyze commutative C*-algebras generated by independent commuting projections (\S~\ref{semigp-comm-proj}, \S~\ref{fam-sub}).

We then consider semigroup C*-algebras and the more general notion of semigroup crossed products by automorphisms (semigroup C*-algebras are the crossed products associated with the trivial action on the complex numbers). We first look at reduced versions (\S~\ref{first-look}). Whenever given a subsemigroup $P$ of a group $G$, there is a canonical $G$-action on a certain C*-algebra associated with every semigroup action of $P$ by automorphisms. We find conditions when the reduced semigroup crossed product by automorphisms embeds as a full corner into the corresponding group crossed product. This leads us to the Toeplitz condition mentioned above. It is introduced and briefly discussed in \S~\ref{sec-T}.

In \S~\ref{var-des}, we then describe reduced and full semigroup crossed products by automorphisms as crossed products by partial automorphisms of inverse semigroups and groupoid crossed products. Here we need to assume that the constructible right ideals of our semigroup are independent. The first main observation is that the Toeplitz condition is precisely what we need to embed full semigroup crossed products by automorphisms as full corners into the corresponding full group crossed products (see Theorem~\ref{thm1}).

As a consequence of our first main result, we determine equivalent characterizations of nuclearity for reduced and full semigroup C*-algebras in \S~\ref{nuc}.

In \S~\ref{ideals}, we study induced ideals of semigroup C*-algebras. We first extend our results on embeddability into full corners and nuclearity to the situation of ideals and quotients (see \S~\ref{ind-ideals}). Induced ideals are obtained from invariant subsets of the spectrum of the canonical commutative subalgebra of the semigroup C*-algebra. Therefore, we explicitly describe this spectrum in \S~\ref{des-spec}. Moreover, we extend the notion of the boundary from \cite{La1} to our general setting. We analyze the boundary action in \S~\ref{bd-act} and find a necessary and sufficient criterion when the boundary quotient is a UCT Kirchberg algebra.

Finally, we turn to examples in \S~\ref{ex}. For quasi-lattice ordered groups, we prove that the analysis from \cite{La-Rae} may be extended to obtain the stronger property of nuclearity of the corresponding semigroup C*-algebras (see \S~\ref{qlo}). We also treat the case of the free product $\Nz_0^{*n}$ in \S~\ref{O}. The boundary quotient in this case is the Cuntz algebra $\cO_n$, and an application of our results yields a description of $\cO_n$ -- up to Morita equivalence -- as a crossed product associated with the action of the free group $\Fz_n$ on the \an{positive part} of its Gromov boundary. Another class of examples is provided by left Ore semigroups (see \S~\ref{Ore}). It turns out that ring C*-algebras are the boundary quotients of the semigroup C*-algebras of the corresponding $ax+b$-semigroups. This explains why several aspects of the structure of ring C*-algebras are very similar to those of the Cuntz algebras.

In \S~\ref{future}, we discuss a few open questions which may be interesting for future research.

I would like to thank M. Laca for bringing \cite{Cr-La1} and \cite{Cr-La2} to my attention. I also thank J. Cuntz for pointing me towards \cite{Kho-Ska1} and \cite{Kho-Ska2}. Moreover, I thank R. Meyer who brought inverse semigroups to my mind.

\section{Preliminaries}

\subsection{The setting}
\label{setting}

Throughout this paper, let $P$ be a subsemigroup of a group $G$. We assume that $P$ contains the unit element $e$ of $G$. All the semigroups in this paper will be unital, and all semigroup homomorphisms shall preserve the units. Moreover, we point out that we are only looking at discrete semigroups and discrete groups.

As explained in \cite{Li2}, the right ideal structure of $P$ plays an important role in the construction and analysis of the semigroup C*-algebras of $P$. By a right ideal of $P$, we mean a subset $X$ of $P$ which is closed under right multiplication, i.e. for all $x \in X$ and $p \in P$, the product $xp$ lies in $X$. Given a subset (for example a right ideal) $X$ of $P$ and a semigroup element $p$, we can form the left translate of $X$ by $p$, i.e. $pX \defeq \menge{px}{x \in X}$, and also the pre-image of $X$ under left multiplication by $p$, i.e. $p^{-1} X \defeq \menge{y \in P}{py \in X}$. Since $G$ also acts on itself by left translations, we can also translate a subset $X$ by a group element $g$. We denote the translation by $g \cdot X \defeq \menge{gx}{x \in X}$. We have for $p$ in $P$ and $X \subseteq P$ that $pX = p \cdot X$, but $p^{-1} X \neq p^{-1} \cdot X$ in general. Instead, we have the relation $p^{-1} X = (p^{-1} \cdot X) \cap P$.

The following family of right ideals was introduced in \cite{Li2}:
\bdefin
Let $\cJ$ be the smallest family of right ideals of $P$ such that
\begin{itemize}
\item $\emptyset, P \in \cJ$;
\item $\cJ$ is closed under left multiplication and pre-images under left multiplication ($X \in \cJ, p \in P \Rarr pX, p^{-1}X \in \cJ$).
\end{itemize}
Elements in $\cJ$ are called constructible right ideals of $P$.
\edefin
As observed in \S~3 in \cite{Li2}, the family $\cJ$ is automatically closed under finite intersections.

In our situation of a subsemigroup of a group, it is also important to consider the following
\bdefin
\label{J_P^G}
Let $\cJ_P^G$ be the smallest family of subsets of $G$ which contains $\cJ$ and which is closed under left translations by group elements ($Y \in \cJ_P^G, g \in G \Rarr g \cdot Y \in \cJ_P^G$) and finite intersections.
\edefin

It is immediate from the definitions that $\cJ$ consists of $\emptyset$ and all right ideals of the form $q_1^{-1} p_1 \dotsm q_n^{-1} p_n P$ ($p_i, q_i \in P$). Moreover, $\cJ_P^G$ is given by all finite intersections of subsets of the form $g \cdot X$, for $g \in G$ and $X \in \cJ$. Actually, $\cJ_P^G$ consists of $\emptyset$ and all finite intersections of subsets of $G$ of the form $g \cdot P$ (for $g \in G$).

\subsection{Semigroups of commuting projections}
\label{semigp-comm-proj}

We will be interested in the following situation: Let $D = C^*(E)$ be a C*-algebra generated by a multiplicative semigroup $E$ of pairwise commuting projections. Given projections $e_1$, ..., $e_n$ in $E$, let $\bigvee_{i=1}^n e_i$ be the smallest projection in $D$ which dominates all the $e_i$ ($1 \leq i \leq n$).

\bdefin
\label{e-ind}
We say that $E$ is independent in $D$ if for all $e$, $e_1$, ..., $e_n$ in $E$, the equation $e = \bigvee_{i=1}^n e_i$ implies that $e=e_i$ for some $1 \leq i \leq n$.
\edefin

Following ideas which appeared in the proof of Proposition~2.24 in \cite{Li2}, we obtain:
\bprop
\label{ind-univ-linind}
The following are equivalent:
\begin{itemize}
\item[(i)] $E$ is independent in $D$,
\item[(ii)] whenever $T$ is a C*-algebra and $\varphi: E \to \Proj(T)$ is a semigroup homomorphism sending $0$ to $0$ if $0 \in E$, then there exists a (unique) homomorphism $D \to T$ given by $e \ma \varphi(e)$ for all $e \in E$,
\item[(iii)] $\menge{e \in E}{e \neq 0}$ is linearly independent in $D$.
\end{itemize}
\eprop
\bproof
We start with \an{(i) $\Rarr$ (ii)}. The idea is to write $D$ as an inductive limit of finite dimensional subalgebras. For every finite subset $F$ of $E$ such that $F \cup \gekl{0}$ is multiplicatively closed, set $D_F = C^*(F) = \lspan(F)$. As all the projections $e \in F$ commute, we may orthogonalize them in $D_F$: For every $0 \neq e \in F$, form the projection $e_{F,D} \defeq e - \bigvee_{e \gneq e' \in F} e'$. As $E$ is independent, all these projections $e_{F,D}$ are non-zero (for $0 \neq e \in F$). Moreover, these projections are pairwise orthogonal, and they generate $D_F$. Thus we obtain $D_F = \bigoplus_{0 \neq e \in F} \Cz \cdot e_{F,D}$. Similarly, form $e_{F,T}^{\varphi} \defeq \varphi(e) - \bigvee_{\varphi(e) \gneq f \in \varphi(F)} f$ in $T$. These projections $e_{F,T}^{\varphi}$ are pairwise orthogonal by construction. Thus there exists by universal property of $\bigoplus_{0 \neq e \in F} \Cz \cdot e_{F,D} \cong \Cz^{\abs{F \setminus \gekl{0}}}$ a homomorphism $D_F \to T$ defined by $e_{F,D} \ma e_{F,T}^{\varphi}$. By construction, $e = \sum_{e \geq e' \in F} e'_{F,D}$ is sent to $\sum_{e \geq e' \in F} {e'}_{F,T}^{\varphi} = \varphi(e)$. Therefore, these homomorphisms $\gekl{D_F \to T}_{F}$ are compatible with the canonical inclusions $D_F \into D_{\ti{F}}$ for $F \subseteq \ti{F}$. Hence they define a homomorphism $D = \overline{\bigcup_F D_F} \to T$ which sends $e \in D$ to $\varphi(e) \in T$ for all $e \in E$, as desired.

For \an{(ii) $\Rarr$ (iii)}, note that by (ii), there exists a homomorphism $D \to D \otimes D, e \ma e \otimes e$. As $D$ is commutative, it does not matter which tensor product we choose. Restricting to $D^{\alg} \defeq \lspan(E)$, we obtain a homomorphism
\bgl
\label{diag-hom}
  D^{\alg} \to D^{\alg} \odot D^{\alg}, e \ma e \otimes e.
\egl
As $D^{\alg}$ is spanned by $E$, we can choose a subset $E'$ of $E$ such that $E'$ is a $\Cz$-basis of $D^{\alg}$. Now take $e \in E$. We can write $e$ as a finite sum $e = \sum_i \lambda_i e'_i$ for some $e'_i \in E'$. The homomorphism from \eqref{diag-hom} sends $e$ to $e \otimes e = \sum_{i,j} \lambda_i \lambda_j e'_i \otimes e'_j$ and $\sum_i \lambda_i e'_i$ to $\sum_i \lambda_i e'_i \otimes e'_i$. But $e$ and $\sum_i \lambda_i e'_i$ coincide, so they have to be sent to the same element. We conclude that
\bgl
\label{sumij}
  \sum_{i,j} \lambda_i \lambda_j e'_i \otimes e'_j = \sum_i \lambda_i e'_i \otimes e'_i.
\egl
As $E'$ is a $\Cz$-basis for $D^{\alg}$, $\menge{e' \otimes e''}{e', e'' \in E'}$ is a $\Cz$-basis for $D^{\alg} \odot D^{\alg}$. Thus we can compare coefficients in \eqref{sumij} and deduce $\lambda_i \lambda_j = 0$ if $i \neq j$ and $\lambda_i^2 = \lambda_i$. It follows that there can at most be one non-zero coefficient $\lambda_i$ which must be $1$. Thus either $e = 0$ or $e = e'_i \in E'$. We deduce that $E' = E \setminus \gekl{0}$. But this means that $\menge{e \in E}{e \neq 0}$ is a $\Cz$-basis of $D^{\alg}$, hence linearly independent.

To see \an{(iii) $\Rarr$ (i)}, we observe that $\bigvee_{i=1}^n e_i = \sum_{\emptyset \neq F \subseteq \gekl{1, \dotsc, n}} (-1)^{\abs{F}-1} \prod_{i \in F} e_i$. Thus an equation of the form $e = \bigvee_{i=1}^n e_i$ with $e \gneq e_i$ for all $1 \leq i \leq n$ would give us a non-trivial relation contradicting linear independence.
\eproof

\subsection{Families of subsets}
\label{fam-sub}

Let us now specialize to a situation which will appear later on in this paper. Let $\fP$ be a discrete set and $\fJ$ be a family of subsets of $\fP$. We assume that $\emptyset \in \fJ$ and that $\fJ$ is closed under finite intersections.

\bdefin
\label{ind}
$\fJ$ is called independent if for all $X$, $X_1, \dotsc, X_n$ in $\fJ$, we have that $X_i \subsetneq X$ for all $1 \leq i \leq n$ implies $\bigcup_{i=1}^n X_i \subsetneq X$.
\edefin

In other words, $\fJ$ is independent if whenever $X = \bigcup_{i=1}^n X_i$ for $X$, $X_1, \dotsc, X_n$ in $\fJ$, then there must be an index $1 \leq i \leq n$ such that $X = X_i$. This independence condition was introduced in \cite{Li2}.

For a subset $X$ of $\fP$, we write $\1z_X$ for the characteristic function of $X$ defined on $\fP$. We view $\1z_X$ as an element of $\ell^{\infty}(\fP)$ and let $\ell^{\infty}(\fP)$ act on $\ell^2(\fP)$ by multiplication operators. Let $E_X$ be the multiplication operator corresponding to $\1z_X$.

\bdefin
We set $D \defeq C^*(\menge{E_X}{X \in \fJ}) \subseteq \ell^{\infty}(\fP) \subseteq \cL(\ell^2(\fP))$.
\edefin
It is easy to see that $\fJ$ is independent if and only if $\menge{E_X}{X \in \fJ}$ is independent in $D$ in the sense of Definition~\ref{e-ind}. Thus, Proposition~\ref{ind-univ-linind} yields in our present setting:
\bcor
\label{ind-equiv}
The following are equivalent:
\begin{itemize}
\item[(i)] $\fJ$ is independent,
\item[(ii)] whenever $T$ is a C*-algebra and $e_X$, $X \in \fJ$, are projections in $T$ satisfying $e_{\emptyset} = 0$ and $e_{X_1 \cap X_2} = e_{X_1} e_{X_2}$ for all $X_1, X_2 \in \fJ$, then there exists a (unique) homomorphism $D \to T$ given by $E_X \ma e_X$ for all $X \in \fJ$,
\item[(iii)] $\menge{E_X}{X \neq \emptyset}$ is linearly independent in $D$.
\end{itemize}
\ecor

From now on, we always assume that $\fJ$ is independent. Let us describe the spectrum of $D$.
\bcor
\label{spec1}
For every function $\phi: \fJ \to \gekl{0,1}$ with $\phi(\emptyset) = 0$ and $\phi(X_1 \cap X_2) = \phi(X_1) \phi(X_2)$ for all $X_1, X_2 \in \fJ$, there exists a unique homomorphism $D \to \Cz$ determined by $E_X \ma \phi(X)$.
\ecor
\bproof
Just set $T = \Cz$ in item (ii) of Corollary~\ref{ind-equiv}.
\eproof

Let us call a subset $\cF$ of $\fJ$ satisfying
\begin{itemize}
\item $X_1 \subseteq X_2 \in \fJ, X_1 \in \cF \Rarr X_2 \in \cF$,
\item $X_1, X_2 \in \cF \Rarr X_1 \cap X_2 \in \cF$,
\item $\emptyset \notin \cF$,
\end{itemize}
a $\fJ$-valued filter.

\bcor
\label{Jvf}
We can identify $\Spec D$ with the set $\Sigma$ of all non-empty $\fJ$-valued filters via $\Spec D \ni \chi \ma \menge{X \in \fJ}{\chi(E_X) = 1}$.
\ecor
\bproof
The inverse of this map is given by sending a non-empty $\fJ$-valued filter to the character $\chi$ of $D$ uniquely determined by $\chi(E_X) = 1$ if $X \in \cF$ and $\chi(E_X) = 0$ if $X \notin \cF$. Such a character exists by Corollary~\ref{spec1}.
\eproof

The topology of pointwise convergence on $\Spec D$ corresponds under the bijection
\bgl
\label{spec2}
  \Spec D \ni \chi \ma \menge{X \in \fJ}{\chi(E_X) = 1} \in \Sigma
\egl
to the following topology on $\Sigma$: For $X$, $X_1, \dotsc, X_n$ in $\fJ$, let
\bgloz
  U(X; X_1, \dotsc, X_n) \defeq \menge{\cF \in \Sigma}{X \in \cF, X_i \notin \cF \fa 1 \leq i \leq n}.
\egloz
Then a basis for the topology on $\Sigma$ induced by the one on $\Spec D$ is given by the open sets
\bgl
\label{basis-top}
  \menge{U(X; X_1, \dotsc, X_n)}{n \in \Zz_{\geq 0}, X, X_1, \dotsc, X_n \in \fJ}.
\egl
Finally, we call a $\fJ$-valued filter which is maximal (in $\Sigma$) with respect to inclusion a $\fJ$-valued ultrafilter.

\bdefin
We let $\Sigma_{\max}$ be the set of all $\fJ$-valued ultrafilters. The subset of $\Spec D$ corresponding to $\Sigma_{\max}$ under the identification \eqref{spec2} is denoted by $(\Spec D)_{\max}$. Moreover, we set $\partial \Sigma \defeq \overline{\Sigma_{\max}} \subseteq \Sigma$ and denote the closed subset of $\Spec D$ corresponding to $\partial \Sigma$ under the homeomorphism~\eqref{spec2} by $\partial \Spec D$.
\edefin

\bremark
\label{F-in-Sigma}
Note that $\cF \in \Sigma$ lies in $\Sigma_{\max}$ if and only if for all $X \in \fJ$, $X \notin \cF$ there is $X' \in \cF$ with $X \cap X' = \emptyset$.
\eremark

\section{A first look at the reduced case}
\label{first-look}

Let us first of all define reduced semigroup C*-algebras and reduced crossed products by automorphisms (see \cite{Li2}). We start with reduced semigroup C*-algebras. Recall that $P$ is a subsemigroup of a group $G$. Let $\menge{\varepsilon_x}{x \in P}$ be the canonical orthonormal basis of $\ell^2(P)$. For every $p \in P$, the formula $V_p \varepsilon_x = \varepsilon_{px}$ extends to an isometry on $\ell^2(P)$. Now the reduced semigroup C*-algebra of $P$ is simply given by the sub-C*-algebra of $\cL(\ell^2(P))$ generated by these isometries $\menge{V_p}{p \in P}$. We denote this concrete C*-algebra by $C^*_r(P)$, i.e. we set
\bdefin
$C^*_r(P) \defeq C^* \rukl{\menge{V_p}{p \in P}} \subseteq \cL(\ell^2(P))$.
\edefin

As we have done in \S~\ref{fam-sub}, we denote by $E_X \in \cL(\ell^2(P))$ the orthogonal projection onto $\ell^2(X) \subseteq \ell^2(P)$ for every subset $X$ of $P$. We then set
\bdefin
$D_r \defeq C^*(\menge{E_X}{X \in \cJ})$.
\edefin
As explained in \cite{Li2}, $D_r$ is a commutative sub-C*-algebra of $C^*_r(P)$.

Now we turn to crossed products. Let $A$ be a C*-algebra which we will always think of as a non-degenerate sub-C*-algebra of $\cL(\cH)$ for some Hilbert space $\cH$. Assume that we are given a $G$-action $\alpha$ on $A$. Define for every $a$ in $A$ the operator $a_{(\alpha \vert_P)} \in \cL(\cH \otimes \ell^2(P))$ by setting $\aalphaP(\xi \otimes \varepsilon_x) = (\alpha_x^{-1} (a) \xi ) \otimes \varepsilon_x$ for all $\xi \in \cH$, $x \in P$.

\bdefin
The reduced automorphic crossed product of $A$ by $P$ is given by $A \rta_{\alpha,r} P \defeq C^*(\menge{\aalphaP (I_{\cH} \otimes V_p)}{a \in A, p \in P}) \subseteq \cL(\cH \otimes \ell^2(P))$ where $I_{\cH}$ is the identity operator on $\cH$.
\edefin

Of course, we can canonically identify $\Cz \rta_{\tr,r} P$ with $C^*_r(P)$.

We now discuss the question whether $A \rta_{\alpha,r} P$ can be embedded as a full corner into an ordinary (reduced) group crossed product. Let $\lambda: G \to \cU(\ell^2(G))$ be the left regular representation of $G$. The group $G$ acts on $\ell^{\infty}(G)$ by left translations. We call this action $\tau$, and we denote the action of $G$ on the multiplication operators on $\ell^2(G)$ corresponding to $\ell^{\infty}(G)$ by $\tau$ as well. It is clear that $\tau$ is spatially implemented by $\lambda$. As before, for a subset $Y$ of $G$, we let $E_Y \in \cL(\ell^2(G))$ be the orthogonal projection onto $\ell^2(Y) \subseteq \ell^2(G)$. In particular, $E_P$ is the orthogonal projection onto $\ell^2(P) \subseteq \ell^2(G)$.

\bdefin
We let $D_P^G$ be the smallest sub-C*-algebra of $\ell^{\infty}(G) \subseteq \cL(\ell^2(G))$ which is $\tau$-invariant and contains $E_P$.
\edefin

\blemma
\label{Dping=span}
With $\cJ_P^G$ from Definition~\ref{J_P^G}, we have $D_P^G = \clspan(\menge{E_Y}{Y \in \cJ_P^G})$.
\elemma
\bproof
Every $Y$ in $\cJ_P^G$ is of the form $\bigcap_{i=1}^n g_i \cdot X_i$ for $g_i \in G$, $X_i \in \cJ$. Thus $E_Y = \prod_{i=1}^n \tau_{g_i}(E_{X_i})$ lies in $D_P^G$. This proves \an{$\supseteq$}. Conversely, the set $\menge{E_Y}{Y \in \cJ_P^G}$ is multiplicatively closed as $\cJ_P^G$ is closed under finite intersections. Moreover, this set is $\tau$-invariant and contains $E_P$. Thus \an{$\subseteq$} holds as well.
\eproof

As in the construction of reduced crossed products, we define for $a \in A$ the operator $a_{(\alpha)} \in \cL(\cH \otimes \ell^2(G))$ by $a_{(\alpha)}(\xi \otimes \varepsilon_x) = (\alpha_x^{-1}(a) \xi) \otimes \varepsilon_x$ for all $\xi$ in $\cH$ and $x$ in $G$. The following is just Proposition~2.5.1 in \cite{C-E-L1} with general coefficients:
\blemma
\label{faithrep}
The homomorphism $A \otimes D_P^G \to \cL(\cH \otimes \ell^2(G))$ determined by $a \otimes d \ma a_{(\alpha)} (I_{\cH} \otimes d)$ and the group homomorphism $G \to \cU(\cH \otimes \ell^2(G)), g \ma I_{\cH} \otimes \lambda_g$ define a covariant representation of $(A \otimes D_P^G, G, \alpha \otimes \tau)$ on $\cH \otimes \ell^2(G)$. The corresponding representation of $(A \otimes D_P^G) \rtimes_{\alpha \otimes \tau,r} G$ is faithful. It sends $(a \otimes d) U_g$ to $a_{(\alpha)} (I_{\cH} \otimes d)(I_{\cH} \otimes \lambda_g)$.
\elemma
Note that since $D_P^G$ is commutative, it does not matter which tensor product $A \otimes D_P^G$ we take.
\bproof
An obvious computation shows that the maps described in the lemma define a covariant representation. Let us show that it gives rise to a faithful representation of the reduced crossed product.

By replacing $\cH$ by $\cH \otimes \ell^2(G)$ and $a \in A$ by $a_{(\alpha)}$, we may without loss of generality assume that the $G$-action $\alpha$ on $A$ is spatially implemented. This means that there exists a group homomorphism $G \to \cU(\cH), g \ma W_g$ such that $\Ad(W_g)(a) = \alpha_g(a)$. We realize the reduced crossed product $(A \otimes D_P^G) \rtimes_{\alpha \otimes \tau,r} G$ as the sub-C*-algebra of $\cL(\cH \otimes \ell^2(G) \otimes \ell^2(G))$ generated by $\menge{(a \otimes d)_{\alpha \otimes \tau} (I_{\cH \otimes \ell^2(G)} \otimes \lambda_g)}{a \in A, d \in D_P^G, g \in G}$ with $(a \otimes d)_{\alpha \otimes \tau} (\xi \otimes \zeta \otimes \varepsilon_x) = (\alpha_x^{-1}(a) \xi) \otimes (\tau_x^{-1}(d) \zeta) \otimes \varepsilon_x$. Now define the unitary
\bgloz
  W: \cH \otimes \ell^2(G) \otimes \ell^2(G) \to \cH \otimes \ell^2(G) \otimes \ell^2(G), 
  \xi \otimes \varepsilon_x \otimes \varepsilon_y \ma W_{x^{-1}} \xi \otimes \varepsilon_{yx} \otimes \varepsilon_{x^{-1}}.
\egloz
A similar computation as in \cite{C-E-L1}, Proposition~2.5.1 shows
\bgloz
  W((a \otimes d)_{\alpha \otimes \tau})W^* = (a_{(\alpha)} (I_{\cH} \otimes d)) \otimes I_{\ell^2(G)}, \ 
  W(I_{\cH \otimes \ell^2(G)} \otimes \lambda_g)W^* = (I_{\cH} \otimes \lambda_g) \otimes I_{\ell^2(G)}.
\egloz
Thus $\Ad(W)$ identifies $(A \otimes D_P^G) \rtimes_{\alpha \otimes \tau,r} G$ with a sub-C*-algebra of $\cL(\cH \otimes \ell^2(G)) \otimes I_{\ell^2(G)}$. Identifying $\cL(\cH \otimes \ell^2(G)) \otimes I_{\ell^2(G)}$ with $\cL(\cH \otimes \ell^2(G))$ in the obvious way, we obtain the desired faithful representation.
\eproof

\bdefin
Let $A \rtimes_{\alpha,r} (\ping)$ be the image of $(A \otimes D_P^G) \rtimes_{\alpha \otimes \tau,r} G$ under the representation from the last lemma. If $A = \Cz$, then we set $C^*_r(\ping) \defeq \Cz \rtimes_{\tr,r} (\ping)$.
\edefin

In the sequel, we denote for $d \in D_P^G$ the canonical multiplier associated with $d$ by $1 \otimes d \in M(A \otimes D_P^G) \subseteq M((A \otimes D_P^G) \rtimes_{\alpha \otimes \tau,r} G$.
\blemma
\label{fullproj}
$1 \otimes E_P$ yields the full corner $(1 \otimes E_P)((A \otimes D_P^G) \rtimes_{\alpha \otimes \tau,r} G)(1 \otimes E_P)$.
\elemma
\bproof
We have to show that $\lspan(((A \otimes D_P^G) \rtimes_{\alpha \otimes \tau,r} G)(1 \otimes E_P)((A \otimes D_P^G) \rtimes_{\alpha \otimes \tau,r} G))$ is dense in $(A \otimes D_P^G) \rtimes_{\alpha \otimes \tau,r} G$.

For every $Y = \bigcap_{i=1}^n g_i \cdot X_i \in \cJ_P^G$ ($g_i \in G$, $X_i \in \cJ$), $a \in A$ and $g \in G$, the operator $(a \otimes E_Y) U_g = (a \otimes E_Y) (\prod_{i=1}^n U_{g_i}(1 \otimes E_{X_i})U_{g_i}^*) U_g = (a \otimes E_Y) (\prod_{i=1}^n U_{g_i}(1 \otimes E_{X_i})(1 \otimes E_P)U_{g_i}^*) U_g$ lies in $((A \otimes D_P^G) \rtimes_{\alpha \otimes \tau,r} G)(1 \otimes E_P)((A \otimes D_P^G) \rtimes_{\alpha \otimes \tau,r} G)$. Here $U_g$ are the canonical unitaries implementing the $G$-action.
\eproof

In the sequel, we do not distinguish between $\cH \otimes \ell^2(P)$ and the subspace $(I_{\cH} \otimes E_P)(\cH \otimes \ell^2(G))$ of $\cH \otimes \ell^2(G)$. In this way, operators on $\cH \otimes \ell^2(P)$ act on $\cH \otimes \ell^2(G)$ (on the orthogonal complement of $(I_{\cH} \otimes E_P)(\cH \otimes \ell^2(G))$, they are simply $0$). For instance, the operator $\aalphaP$ is the same as $(I_{\cH} \otimes E_P) a_{(\alpha)} (I_{\cH} \otimes E_P)$ and $I_{\cH} \otimes V_p$ is nothing else but $(I_{\cH} \otimes E_P)(I_{\cH} \otimes \lambda_p)(I_{\cH} \otimes E_P)$ for all $p \in P$. As $A \rta_{\alpha,r} P$ is the C*-algebra generated by $\aalphaP (I_{\cH} \otimes V_p)$ ($a \in A$, $p \in P$), we see that $A \rta_{\alpha,r} P$ is (or can be, in the way explained above, canonically identified with) a sub-C*-algebra of the full corner $(I_{\cH} \otimes E_P)(A \rtimes_{\alpha,r}(\ping))(I_{\cH} \otimes E_P)$. We now address the question when these two C*-algebras are actually the same.

\blemma
\label{AP=APG}
The following statements are equivalent:
\begin{enumerate}
\item[(i)] We have $A \rta_{\alpha,r} P = (I_{\cH} \otimes E_P)(A \rtimes_{\alpha,r}(\ping))(I_{\cH} \otimes E_P)$ for every C*-dynamical system $(A,G,\alpha)$,
\item[(ii)] $C^*_r(P) = E_P C^*_r(\ping) E_P$,
\item[(iii)] for all $g \in G$, $E_P \lambda_g E_P$ lies in $C^*_r(P)$;
\end{enumerate}
and either of these statements implies
\begin{itemize}
\item[(iv)] $D_r \supseteq E_P D_P^G E_P$.
\end{itemize}
\elemma
\bproof
\an{(i) $\Rarr$ (ii)} is trivial.

\an{(ii) $\Rarr$ (iii)}: $E_P \lambda_g E_P = E_P (E_P \lambda_g) E_P \in E_P C^*_r(\ping) E_P = C^*_r(P)$.

\an{(iii) $\Rarr$ (iv)}: By Lemma~\ref{Dping=span} and the definition of $\cJ_P^G$ from \S~\ref{setting}, it suffices to prove that $E_P E_{g \cdot X} E_P$ lies in $D_r$ for all $g \in G$ and $X \in \cJ$. First of all, $E_P E_{g \cdot X} E_P = E_P \lambda_g E_X \lambda_g^* E_P = (E_P \lambda_g E_P) E_X (E_P \lambda_g E_P)^*$ lies in $C^*_r(P)$ by (iii). Moreover, $E_P E_{g \cdot X} E_P$ is obviously contained in $\ell^{\infty}(P)$ viewed as multiplication operators on $\ell^2(P)$. Thus $E_P E_{g \cdot X} E_P$ lies in $C^*_r(P) \cap \ell^{\infty}(P)$, and $C^*_r(P) \cap \ell^{\infty}(P) = D_r$ by Remark~3.12 in \cite{Li2}.

\an{(iii) \& (iv) $\Rarr$ (i)}: We have to show that for every $a \in A$, $Y \in \cJ_P^G$ and $g \in G$, $(I_{\cH} \otimes E_P)(a_{(\alpha)} (I_{\cH} \otimes E_Y)(I_{\cH} \otimes \lambda_g))(I_{\cH} \otimes E_P)$ lies in $A \rta_{\alpha,r} P$. We have
\bglnoz
  && (I_{\cH} \otimes E_P)(a_{(\alpha)} (I_{\cH} \otimes E_Y)(I_{\cH} \otimes \lambda_g))(I_{\cH} \otimes E_P) \\
  &=& (I_{\cH} \otimes E_P)a_{(\alpha)}(I_{\cH} \otimes E_P) (I_{\cH} \otimes E_P E_Y \lambda_g E_P) \\
  &=& \aalphaP (I_{\cH} \otimes (E_P E_Y E_P) (E_P \lambda_g E_P)).
\eglnoz
But $E_P E_Y E_P$ lies in $D_r$ by (iv) and $E_P \lambda_g E_P$ is in $C^*_r(P)$ by (iii). Since $\aalphaP(I_{\cH} \otimes C^*_r(P))$ lies in $A \rta_{\alpha,r} P$, we are done.
\eproof

Let us now summarize what we have obtained so far. Combining Lemmas~\ref{faithrep}, \ref{fullproj} and \ref{AP=APG}, we obtain
\bcor
\label{full-red}
If $\ping$ satisfies one of the equivalent conditions (i), (ii) or (iii) from Lemma~\ref{AP=APG}, then the homomorphism $A \rta_{\alpha,r} P \to (A \otimes D_P^G) \rtimes_{\alpha \otimes \tau,r} G$ determined by $\aalphaP (I_{\cH} \otimes V_p) \ma (a \otimes E_P) (1 \otimes E_P) U_p (1 \otimes E_P)$ identifies $A \rta_{\alpha,r}P$ with the full corner $(1 \otimes E_P)((A \otimes D_P^G) \rtimes_{\alpha \otimes \tau,r} G)(1 \otimes E_P)$ of $(A \otimes D_P^G) \rtimes_{\alpha \otimes \tau,r} G$.
\ecor

\section{The Toeplitz condition}
\label{sec-T}

We now introduce a condition on the inclusion $\ping$ which is (at least a priori) stronger then (iii) from Lemma~\ref{AP=APG}.
\bdefin
\label{Toeplitz}
We say that $\ping$ satisfies the Toeplitz condition (or simply that $\ping$ is Toeplitz) if for every $g \in G$ with $E_P \lambda_g E_P \neq 0$, there exist $p_1, q_1, \dotsc, p_n, q_n$ in $P$ such that $E_P \lambda_g E_P = V_{p_1}^* V_{q_1} \dotsm V_{p_n}^* V_{q_n}$.
\edefin
This Toeplitz condition will play an important role in the next section, when we consider full versions. Moreover, we will see that in examples, the Toeplitz condition will naturally appear. In addition, it has the following consequences:
\blemma
\label{lemma-T}
If $\ping$ is Toeplitz, then
\begin{enumerate}
\item[(i)] for all $g$ in $G$ and $X$ in $\cJ$, $P \cap (g \cdot X)$ lies in $\cJ$,
\item[(ii)] $\cJ_P^G = \menge{g \cdot X}{g \in G, X \in \cJ}$ (i.e. intersections are not needed).
\end{enumerate}
If $\cJ$ is independent and $\ping$ is Toeplitz, then
\begin{itemize} 
\item[(iii)] $\cJ_P^G$ is independent.
\end{itemize}
\elemma
\bproof
Given $g \in G$ with $E_P \lambda_g E_P \neq 0$, the Toeplitz condition says that there exist $p_1$, $q_1$, ..., $p_n$, $q_n$ in $P$ such that $E_P \lambda_g E_P = V_{p_1}^* V_{q_1} \dotsm V_{p_n}^* V_{q_n}$. This implies that
\bglnoz
  E_{P \cap (g \cdot X)} &=& (E_P \lambda_g E_P) E_X (E_P \lambda_g E_P)^* = V_{p_1}^* V_{q_1} \dotsm V_{p_n}^* V_{q_n} E_X V_{q_n}^* V_{p_n} \dotsm V_{q_1}^* V_{p_1} \\
  &=& \E{p_1^{-1} q_1 \dotsm p_n^{-1} q_n X}.
\eglnoz
Thus we deduce $P \cap (g \cdot X) = p_1^{-1} q_1 \dotsm p_n^{-1} q_n X \in \cJ$. If $E_P \lambda_g E_P = 0$, then the computation shows that $P \cap (g \cdot X) = \emptyset$ lies in $\cJ$. This proves (i). To prove (ii), we just have to show that the right hand side in (ii) is closed under finite intersections. Take $g_1$, $g_2$ in $G$ and $X_1$, $X_2$ in $\cJ$. Then $(g_1 \cdot X_1) \cap (g_2 \cdot X_2) = g_1 \cdot (X_1 \cap ((g_1^{-1} g_2) \cdot X_2)) = g_1 \cdot (X_1 \cap \underbrace{P \cap (g_1^{-1} g_2) \cdot X_2}_{\in \cJ \text{ by (i)}})$ is of the desired form by (i).

Now let us prove (iii) assuming that $\cJ$ is independent and that $\ping$ is Toeplitz. By (ii), it suffices to prove that given $g$, $g_1, \dotsc, g_n$ in $G$ and $X$, $X_1, \dotsc X_n$ in $\cJ$ such that $g \cdot X = \bigcup_{i=1}^n g_i \cdot X_i$, we must have $g \cdot X = g_i \cdot X_i$ for some $i$. Now $g \cdot X = \bigcup_{i=1}^n g_i \cdot X_i$ implies $X = \bigcup_{i=1}^n (g^{-1} g_i) \cdot X_i$. In particular, since $X \subseteq P$, we must have $(g^{-1} g_i) \cdot X_i \subseteq P$ for all $1 \leq i \leq n$. Therefore $(g^{-1} g_i) \cdot X_i = P \cap ((g^{-1} g_i) \cdot X_i)$ lies in $\cJ$ by (i). As $\cJ$ is independent, there exists $i$ such that $X = (g^{-1} g_i) \cdot X_i$. Thus $g \cdot X = g_i \cdot X_i$.
\eproof

\section{Various descriptions of semigroup crossed products by automorphisms}
\label{var-des}

\subsection{The full versions}

We now turn to full semigroup C*-algebras and full crossed products by automorphisms. We work with the version of full semigroup C*-algebras from \cite{Li2}, \S~3. Recall that $P$ is a subsemigroup of the group $G$.
\bdefin
\label{full-semigpC}
Let $C^*_s(P)$ be the universal C*-algebra generated by isometries $\menge{v_p}{p \in P}$ and projections $\menge{e_X}{X \in \cJ}$ satisfying the following relations:
\begin{enumerate}
\item[I.] $v_p v_q = v_{pq}$ for all $p$, $q$ in $P$,
\item[II.] $e_{\emptyset} = 0$,
\item[III.] whenever $p_1, q_1, \dotsc, p_n, q_n \in P$ satisfy $p_1^{-1} q_1 \dotsm p_n^{-1} q_n = e$ in $G$, then $$v_{p_1}^* v_{q_1} \dotsm v_{p_n}^* v_{q_n} = \e{q_n^{-1} p_n \dotsm q_1^{-1} p_1 P}.$$
\end{enumerate}
We set $D \defeq C^*(\menge{e_X}{X \in \cJ}) \subseteq C^*_s(P)$.
\edefin
These relations are satisfied in $C^*_r(P)$ by \cite{Li2}, Lemma~3.1. Therefore we obtain a homomorphism, the left regular representation, $\lambda: C^*_s(P) \to C^*_r(P)$ given by $v_p \ma V_p$ and $e_X \ma E_X$. As observed in \cite{Li2}, Lemma~3.3., the map $\cJ \ni X \ma e_X \in C^*_s(P)$ is a spectral measure, i.e. $e_P = 1$ and $e_{X_1 \cap X_2} = e_{X_1} e_{X_2}$.

We now define full crossed products by automorphisms. Let $A$ be a C*-algebra and $\alpha$ a $G$-action on $A$. In a certain sense, we now form a universal model of the reduced crossed product $A \rta_{\alpha,r} P$.
\bdefin
The full crossed product of $(A,P,\alpha)$ is a C*-algebra $A \rta_{\alpha,s} P$ which comes with two homomorphisms $\iota: A \to A \rta_{\alpha,s} P$ and $\overline{(\cdot)}: C^*_s(P) \to M(A \rta_{\alpha,s} P)$, $x \ma \overline{x}$, with $\iota(\alpha_p(a)) \overline{v_p} = \overline{v_p} \iota(a)$ for all $p \in P$, $a \in A$, such that the following universal property holds:

Whenever $T$ is a C*-algebra, $\iota': A \to T$ and $(\cdot)': C^*_s(P) \to M(T)$, $x \ma x'$, are homomorphisms satisfying 
\bgl
\label{cov-rel}
\iota'(\alpha_p(a)) v_p' = v_p' \iota'(a) \fa p \in P, \: a \in A,
\egl
then there exists a unique homomorphism $\iota' \rtimes (\cdot)': A \rta_{\alpha,s} P \to T$ sending $\iota(a) \overline{x}$ to $\iota'(a) x'$ for all $a \in A$ and $x \in C^*_s(P)$.
\edefin
Existence of $(A \rta_{\alpha,s} P, \iota, \overline{(\cdot)})$ follows from the existence of Murphy's crossed product (see \cite{Mur2}, \S~1) and the observation that our construction is -- in a canonical way -- a quotient of Murphy's. Moreover, it is clear that $(A \rta_{\alpha,s} P, \iota, \overline{(\cdot)})$ is unique up to canonical isomorphism.

The homomorphisms $A \to A \rta_{\alpha,r} P$, $a \ma \aalphaP$ and $C^*_s(P) \to M(A \rta_{\alpha,r} P)$, $x \ma I_{\cH} \otimes \lambda(x)$ satisfy the covariance relation \eqref{cov-rel}. Thus, by universal property of $A \rta_{\alpha,s} P$, there exists a homomorphism $\lambda_{(A,P,\alpha)}: A \rta_{\alpha,s} P \to A \rta_{\alpha,r} P$ sending $\iota(a) \overline{x}$ to $\aalphaP (I_{\cH} \otimes \lambda(x))$.

Of course, in case $A = \Cz$ we can canonically identify $\Cz \rta_{\tr,s} P$ with $C^*_s(P)$ so that $\lambda_{(\Cz,P,\tr)}$ becomes the left regular representation $\lambda$.

\subsection{Inverse semigroups of partial isometries}
\label{invsemi}

\bdefin
Let $S$ be the multiplicative subsemigroup of $C^*_s(P)$ generated by the isometries $v_q$ and their adjoints $v_p^*$, i.e.
\bgloz
  S \defeq \menge{v_{p_1}^* v_{q_1} \dotsm v_{p_n}^* v_{q_n}}{n \in \Zz_{\geq 0}; p_i, q_i \in P} \cup \gekl{0} \subseteq C^*_s(P).
\egloz
Also, in the reduced case, let $S_r$ be the corresponding subsemigroup of $C^*_r(P)$, i.e. 
\bgloz
  S_r \defeq \menge{V_{p_1}^* V_{q_1} \dotsm V_{p_n}^* V_{q_n}}{n \in \Zz_{\geq 0}; p_i, q_i \in P} \cup \gekl{0} = \lambda(S) \subseteq C^*_r(P).
\egloz
\edefin
It is clear that $S$ and $S_r$ are *-invariant semigroups of partial isometries with commuting range and source projections, hence inverse semigroups.

\blemma
\label{g}
The map $g_r: S_r \setminus \gekl{0} \to G, V_{p_1}^* V_{q_1} \dotsm V_{p_n}^* V_{q_n} \ma p_1^{-1} q_1 \dotsm p_n^{-1} q_n$ is well-defined. For $0 \neq V \in S_r$, $g_r(V)$ is determined by the property that for every $x \in P$, $V \ve_x \neq 0 \Rarr V \ve_x = \ve_{g_r(V)x}$. Moreover, we have $g_r(V^*) = (g_r(V))^{-1}$ for $0 \neq V \in S_r$ and $g_r(V_1 V_2) = g_r(V_1) g_r(V_2)$ for $V_1$, $V_2$ in $S_r$ such that $V_1 V_2 \neq 0$.
\elemma
\bproof
For every $0 \neq V \in S_r$, we obviously have for every $x \in P$ that $V \ve_x$ is either $0$ or of the form $\ve_{p_1^{-1} q_1 \dotsm p_n^{-1} q_n x}$ if $V = V_{p_1}^* V_{q_1} \dotsm V_{p_n}^* V_{q_n}$.
\eproof

This lemma allows the following
\bdefin
We set $g \defeq g_r \circ \lambda: S \setminus \gekl{0} \to G, v_{p_1}^* v_{q_1} \dotsm v_{p_n}^* v_{q_n} \ma p_1^{-1} q_1 \dotsm p_n^{-1} q_n$.
\edefin

\blemma
\label{lambda-inj}
If $\cJ$ is independent, then $\lambda: S \to C^*_r(P), s \ma \lambda(s)$ is injective, or in other words, $\lambda$ identifies $S$ with $S_r$.
\elemma
\bproof
Take two elements $s_1$, $s_2$ from $S$ with $s_1 \neq s_2$, and assume without loss of generality $s_1 \neq 0$, hence $s_1^* s_1 \neq 0$. As $g(s_1^* s_1) = e$, $s_1^* s_1$ lies in $D$ by relation III in Definition~\ref{full-semigpC}. As $\cJ$ is independent, $\lambda$ is injective on $D$ by \cite{Li2}, Corollary~3.4. Thus $\lambda(s_1^* s_1) \neq 0$, hence also $\lambda(s_1) \neq 0$. So if $s_2 = 0$, we conclude that $\lambda(s_1) \neq \lambda(s_2)$. We may now assume that $s_1 \neq 0$ and $s_2 \neq 0$. We start with the case $g(s_1) \neq g(s_2)$. There exists $x \in P$ such that $\lambda(s_1) \ve_x = \ve_{g(s_1)x}$. As $\lambda(s_2) \ve_x$ is either $0$ or equal to $\ve_{g(s_2)x} \neq \ve_{g(s_1)x}$, we have $\lambda(s_1) \neq \lambda(s_2)$. If $g(s_1) = g(s_2)$, then $(s_1 - s_2)^* (s_1 - s_2)$ lies in $D$ by relation III in Definition~\ref{full-semigpC}. As $\lambda$ is injective on $D$ ($\cJ$ is assumed to be independent) and $s_1 \neq s_2$ by assumption, we conclude that $\lambda((s_1 - s_2)^* (s_1 - s_2)) \neq 0$, hence $\lambda(s_1 - s_2) \neq 0$.
\eproof

\bremark
We mention that both $S$ and $S_r$ can be identified (up to $0$) with the left inverse hull of $P$ (see \cite{Nor}). This gives a purely algebraic description of these inverse semigroups in terms of $P$.
\eremark

\subsection{Crossed products by partial automorphisms}
\label{cropro-partialauto}

Our goal is to describe $A \rta_{\alpha,s} P$ as a crossed product of $A \otimes D$ by $S$. The reader may consult \cite{Sie} for the general construction of crossed products associated with partial actions of inverse semigroups.

First of all, it is easy to see that for every $s \in S$, we have an automorphism
\bgloz
  \beta_s: A \otimes s^*Ds \to A \otimes sDs^*, \: a \otimes d \ma \alpha_{g(s)}(a) \otimes sds^*.
\egloz
(For $s=0$, we let $\beta_0$ be the zero map $\gekl{0} \to \gekl{0}$.) Moreover, we have $s^*Ds = s^*ss^*Ds \subseteq s^*sD = s^*sDs^*s \subseteq s^*Ds$ so that $s^*Ds = s^*sD$ is an ideal of $D$. In this way, $S$ acts on $A \otimes D$ by partial automorphisms, i.e. we have a semigroup homomorphism $S \to \PAut(A \otimes D)$, $s \ma \beta_s$.

\bprop
\label{cd1}
We can identify $A \rta_{\alpha,s} P$ and $(A \otimes D) \rtimes_{\beta} S$ by mutually inverse homomorphisms
\bglnoz
  && A \rta_{\alpha,s} P \to (A \otimes D) \rtimes_{\beta} S, \: \iota(a) \overline{s} \ma (a \otimes ss^*) \delta_s, \\
  && (A \otimes D) \rtimes_{\beta} S \to A \rta_{\alpha,s} P, \: (a \otimes ss^*) \delta_s \ma \iota(a) \overline{s}.
\eglnoz
\eprop
\bproof
We use the universal properties of these two crossed products to show existence of these homomorphisms. To construct the homomorphism $A \rta_{\alpha,s} P \to (A \otimes D) \rtimes_{\beta} S$, represent $(A \otimes D) \rtimes_{\beta} S$ faithfully and non-degenerately on a Hilbert space $H$. Proposition~4.7 in \cite{Sie} yields representations of $A \otimes D$ and $S$ on $H$ such that \eqref{cov-rel} is satisfied. By universal property of $A \rta_{\alpha,s} P$, these representations give rise to the desired homomorphism $A \rta_{\alpha,s} P \to (A \otimes D) \rtimes_{\beta} S$, $\iota(a) \overline{s} \ma (a \otimes ss^*) \delta_s$.

In the reverse direction, the homomorphisms $A \otimes D \to A \rta_{\alpha,s} P$, $a \otimes d \ma \iota(a) \overline{d}$ and $S \to M(A \rtaas P)$, $s \ma \overline{s}$, form a covariant representation of $(A \otimes D,S,\beta)$ in the sense of \cite{Sie}, Definition~3.4 (having represented $A \rtaas P$ faithfully and non-degenerately on a Hilbert space). Then the universal property of $(A \otimes D) \rtimes_{\beta} S$ (see \cite{Sie}, Proposition~4.8) gives the desired homomorphism $(A \otimes D) \rtimes_{\beta} S \to A \rta_{\alpha,s} P$ sending $(a \otimes ss^*) \delta_s$ to $\iota(a) \overline{s}$.

Finally, it is immediate that these homomorphisms are mutual inverses.
\eproof

We can also consider reduced versions. Let us first define the left regular representation of $S$. Set $S\reg = S \setminus \gekl{0}$ and let $\menge{\ve_x}{x \in S\reg}$ be the canonical orthonormal basis of $\ell^2(S\reg)$. Define for every $s \in S$ the partial isometry $\lambda_S(s)$ on $\ell^2(S\reg)$ by
\bgloz
  \lambda_S(s) \ve_x = 
  \bfa
    \ve_{sx} \falls x = s^* s x \\
    0 \sonst.
  \efa
\egloz
Moreover, for every $a \in A$, let $a_{\alpha,S}$ be the operator on $\cH \otimes \ell^2(S\reg)$ given by $a_{\alpha,S}(\xi \otimes \ve_x) = \alpha_{g(x)}^{-1}(a) \xi \otimes \ve_x$. The homomorphisms $A \otimes D \ni a \otimes s^* s \ma a_{\alpha,S} \cdot (I_{\cH} \otimes \lambda_S(s^* s)) \in \cL(\cH \otimes \ell^2(S\reg))$ and $S \ni s \ma I_{\cH} \otimes \lambda_S(s) \in \cL(\cH \otimes \ell^2(S\reg))$ form a covariant representation of $(A \otimes D,S,\beta)$ in the sense of Definition~3.4 in \cite{Sie}. So we obtain a representation of $(A \otimes D) \rtimes_{\beta} S$ on $\cH \otimes \ell^2(S\reg)$. Its image is the sub-C*-algebra of $\cL(\cH \otimes \ell^2(S\reg))$ generated by $a_{\alpha,S} \cdot (I_{\cH} \otimes \lambda_S(s))$ ($a \in A$, $s \in S$), and we denote this C*-algebra by $(A \otimes D) \rtimes_{\beta,r} S$:
\bdefin
$(A \otimes D) \rtimes_{\beta,r} S \defeq C^*(\menge{a_{\alpha,S} \cdot (I_{\cH} \otimes \lambda_S(s))}{a \in A, s \in S}) \subseteq \cL(\cH \otimes \ell^2(S\reg))$.
\edefin
\blemma
\label{reduced-invsemigp-cropro}
There is a canonical homomorphism $(A \otimes D) \rtimes_{\beta,r} S \to A \rta_{\alpha,r} P$ sending $a_{\alpha,S} \cdot (I_{\cH} \otimes \lambda_S(s))$ to $\aalphaP \lambda(s)$. If $\cJ$ is independent, then this homomorphism is an isomorphism.
\elemma
\bproof
This is just the analogue of Corollary~3.2.13 and Theorem~3.2.14 from \cite{Nor} for general coefficients. The same proof as in \cite{Nor} works here as well.
\eproof

\subsection{Groupoids associated with inverse semigroups}
\label{groupoids}

To every inverse semigroup belongs a groupoid. The reader may consult \cite{Pat}, \S~4.3 or \cite{Kho-Ska1} for the general construction. Note that all the groupoids in this paper will be $r$-discrete (also called \'{e}tale) and Hausdorff. We assume in the rest of this section (\S~\ref{var-des}) that $\cJ$ is independent (see Definition~\ref{ind}) and that $\ping$ satisfies the Toeplitz condition (see Definition~\ref{Toeplitz}).

Let us now explain how to construct the groupoid for our specific inverse semigroup $S$. Set $E \defeq \menge{s^*s}{s \in S}$, this is the set of idempotents of $S$. In our case, $E = \menge{e_X}{X \in \cJ}$. The unit space of the groupoid of $S$ is given by the semicharacters on $E$. This unit space can be canonically identified with $\Spec D$. As $\cJ$ is independent, we can identify $D$ and $D_r$ using Corollary~3.4 from \cite{Li2}. Thus $\lambda$ induces an identification of $\Spec D$ with $\Omega \defeq \Spec(D_r)$.

To form the groupoid of $S$, we take $\hat{S} \defeq \menge{(s,\chi) \in S \times \Omega}{\chi(\lambda(s^*s))=1}$ equipped with the subspace topology of $S \times \Omega$. Here $S$ is viewed as a discrete set, and $\Omega$ carries the usual topology of pointwise convergence. Next we define an equivalence relation $\sim$ on $\hat{S}$ by setting
\bgloz
  (s_1,\chi_1) \sim (s_2,\chi_2) :\LRarr \chi_1 = \chi_2 \text{ and there is } e \in E \text{ with } \chi_1(\lambda(e)) = 1 \text{ and } s_1 e = s_2 e.
\egloz
Then the groupoid of $S$ is defined by $\cG(S) \defeq \hat{S} / \sim$ with the quotient topology induced from $\hat{S}$. We write $[s,\chi]$ for the equivalence class of $(s,\chi) \in \hat{S}$. The groupoid structure of $\cG(S)$ is given as follows: First, for $(s,\chi) \in \hat{S}$, let $s.\chi$ be the character $\chi(\lambda(s^* \sqcup s)) = \chi \circ \Ad(\lambda(s)^*)$. Two elements $[s_1,\chi_1]$ and $[s_2,\chi_2]$ of $\cG(S)$ are composable if $s_2.\chi_2 = \chi_1$, and in that case, the product is given by $[s_1,\chi_1][s_2,\chi_2] = [s_1s_2,\chi_2]$. The inverse map is given by $[s,\chi]^{-1} = [s^*,s.\chi]$. Moreover, the range and source maps are $r: \cG(S) \to \Omega$, $[s,\chi] \ma s.\chi$ and $d: \cG(S) \to \Omega$, $[s,\chi] \ma \chi$.

Let us now compare this groupoid with another one. Namely, we have a canonical transformation groupoid associated with the dynamical system $(D_P^G,G,\tau)$ since $D_P^G$ is commutative. Set $\Omega_P^G \defeq \Spec(D_P^G)$. The group $G$ acts on $\Omega_P^G$ from the right by $\chi g = \chi \circ \tau_g$ (this is just the transpose of $\tau$). The corresponding transformation groupoid is denoted by $\cG \defeq \Omega_P^G \rtimes G$. As a topological space, $\cG$ is simply the product space $\Omega_P^G \times G$. Two elements $(\chi_1,g_1)$ and $(\chi_2,g_2)$ are composable if $\chi_1 g_1 = \chi_2$, and in this case we have $(\chi_1,g_1)(\chi_2,g_2) = (\chi_1,g_1g_2)$. The inverse map is given by $(\chi,g)^{-1} = (\chi g,g^{-1})$. Furthermore, the range and source maps are given by $r: \cG \to \Omega_P^G$, $(\chi,g) \ma \chi$ and $d: \cG \to \Omega_P^G$, $(\chi,g) \ma \chi g$.

Now we restrict $\cG$ to a subset of $\Omega_P^G$. By our assumption that $\ping$ is Toeplitz, we know that $D_r = E_P D_P^G E_P$ by Lemma~\ref{AP=APG}. Therefore we can define a surjective homomorphism $c: D_P^G \to D_r$, $x \ma E_P x E_P$. This homomorphism induces an embedding $c^*: \Omega \to \Omega_P^G$, $\chi \ma \chi \circ c$. We set $N \defeq c^*(\Omega)$.

\blemma
\label{N}
We have $N = \menge{\chi \in \Omega_P^G}{\chi(E_P) = 1}$.
\elemma
\bproof
\an{$\subseteq$}: $(\chi \circ c) (E_P) = \chi(E_P) = 1$ for all $\chi \in \Omega$ as $E_P$ is the unit of $D_r$.

\an{$\supseteq$}: If $\chi \in \Omega_P^G$ satisfies $\chi(E_P) = 1$, then $\chi(x) = \chi(E_P) \chi(x) \chi(E_P) = (\chi \circ c) (x)$ for all $x \in D_P^G$. Thus $\chi = (\chi \vert_{D_r}) \circ c \in N$.
\eproof

\bcor
\label{c-open}
$N$ is clopen in $\Omega_P^G$ and $c^*: \Omega \to \Omega_P^G$ is open.
\ecor
\bproof
The first assertion is immediate from the previous lemma. To see our second claim, observe that $c^* \vert^N: \Omega \to N$ is a homeomorphism, being a continuous bijection between compact Hausdorff spaces. Now given an open subset $U \subseteq \Omega$, $c^*(U)$ is open in $N$, hence also in $\Omega_P^G$ as $N$ is open.
\eproof

We now form the groupoid
\bgl
\label{GNN}
  \cG_N^N \defeq r^{-1}(N) \cap d^{-1}(N).
\egl
$\cG_N^N$ inherits from $\cG$ the structure of a topological groupoid by taking the subspace topology and restricting the product and the inverse map.

The next observation tells us that restricting to $N$ does not lead so far away:
\blemma
$N$ meets every orbit in $\cG^{(0)}$, i.e. for every $\chi \in \cG^{(0)} = \Omega_P^G$, there exists $g \in G$ such that $d(\chi,g)$ lies in $N$. Moreover, the restricted range and source maps $r \vert_{d^{-1}(N)}: d^{-1}(N) \to \Omega_P^G$ and $d \vert_{d^{-1}(N)}: d^{-1}(N) \to N$ are open.
\elemma
\bproof
For every $\chi \in \Omega_P^G$ there exists $Y \in \cJ_P^G$ such that $\chi(E_Y) = 1$. This subset $Y$ of $G$ must be of the form $Y = \bigcap_{i=1}^n g_i \cdot X_i$ for some $n \geq 1$, $g_i \in G$ and $X_i \in \cJ$. Thus $\chi(E_Y) = 1$ implies $\chi(E_{g_1 \cdot X_1})=1$, hence $(\chi g_1)(E_{X_1})=1$. As $E_{X_1} \leq E_P$, we conclude that $(\chi g_1)(E_P)=1$, which means that $d(\chi,g_1) = \chi g_1$ lies in $N$ by Lemma~\ref{N}.

To see that the restricted range and source maps are open, take an open subset $U$ of $\cG$. Then $r(U \cap d^{-1}(N))$ is open as $U$ and $d^{-1}(N)$ are open (recall that $N$ is open) and $r$ is an open map from $\cG$ to $\cG^{(0)}$. Also, $d(U \cap d^{-1}(N)) = d(U) \cap N$ is open in $N$ as $U$ is open and $d$ is an open map $\cG \to \cG^{(0)}$.
\eproof

Setting $\cG_N \defeq d^{-1}(N)$, we have
\bcor
\label{GGNN}
$\cG_N$ together with the restricted range and source maps and the left $\cG$-action and the right $\cG_N^N$-action induced by the product in $\cG$ is a $(\cG,\cG_N^N)$-equivalence in the sense of \cite{M-R-W}.
\ecor
\bproof
This follows from the previous lemma using Example~2.7 in \cite{M-R-W}.
\eproof

Now we return to the groupoid $\cG(S)$ and compare it with $\cG_N^N$.
\bprop
\label{id-groupoids}
Under our standing assumptions that $\cJ$ is independent and $\ping$ is Toeplitz, we can identify $\cG(S)$ with $\cG_N^N$ as topological groupoids via
\bgloz
  \Phi: \cG(S) \to \cG_N^N, \: [s,\chi] \ma ((c^* \chi) g(s)^{-1},g(s)).
\egloz
\eprop
\bproof
First of all, $\Phi$ is well-defined: Namely, $(s_1,\chi) \sim (s_2,\chi)$ implies that there exists $X \in \cJ$ such that $\chi(E_X)=1$ and $s_1 e_X = s_2 e_X$. Thus $\chi(\lambda(e_X s_1^* s_1 e_X))=1$, and we conclude that $s_1 e_X = s_2 e_X$ is not zero. Using Lemma~\ref{g}, we see that $g(s_1) = g(s_1 e_X) = g(s_2 e_X) = g(s_2)$. Therefore, $((c^* \chi) g(s)^{-1},g(s))$ really only depends on the equivalence class of $(s,\chi) \in \hat{S}$.

To see that $((c^* \chi) g(s)^{-1},g(s))$ lies in $\cG_N^N$, we have to check that the range and source of $((c^* \chi) g(s)^{-1},g(s))$ lie in $N$. For the source this is obvious. To see that $r((c^* \chi) g(s)^{-1},g(s)) = (c^* \chi) g(s)^{-1}$ lies in $N$, we show
\bgl
\label{cschi}
  c^*(s.\chi) = (c^* \chi)g(s)^{-1} \fa (s,\chi) \in \hat{S}.
\egl
Write $s = v_{p_1}^* v_{q_1} \dotsm v_{p_n}^* v_{q_n}$. For $Y$ in $\cJ_P^G$, we compute that $\lambda(s^*) E_P E_Y E_P \lambda(s) = E_{q_n^{-1} p_n \dotsm q_1^{-1} p_1 (Y \cap P)} = E_{(g(s)^{-1} \cdot Y) \cap P} E_{q_n^{-1} p_n \dotsm q_1^{-1} p_1 P}$. Thus we deduce $c^*(s.\chi)(E_Y) = \chi(\lambda(s^*) E_P E_Y E_P \lambda(s)) = \chi(E_P \tau_{g(s)^{-1}}(E_Y) E_P) \chi(\lambda(s^*s)) = ((c^* \chi)g(s)^{-1})(E_Y)$. This proves \eqref{cschi}.

It is clear that \eqref{cschi} implies $(c^* \chi) g(s)^{-1} \in N$. So far, we have shown that $\Phi$ is well-defined.

To show that $\Phi$ is injective, take $[s_1,\chi_1]$ and $[s_2,\chi_2]$ from $\cG(S)$. Assume that $\Phi([s_1,\chi_1]) = \Phi([s_2,\chi_2]) = (\chi,g)$. Then we must have $g(s_1)=g(s_2)=g$. Moreover, $c^* \chi_1$ and $c^* \chi_2$ must coincide with $\chi g$. The equality $c^* \chi_1 = c^* \chi_2$ implies $\chi_1 = \chi_2$ as $c^*$ is injective. Finally, to prove $[s_1,\chi_1] = [s_2,\chi_2]$, we observe that $(\chi g)(\lambda(s_1^* s_1)) = (\chi g)(\lambda(s_2^* s_2)) = 1$ implies $(\chi g)(\lambda(s_1^* s_1 s_2^* s_2)) = 1$. Now set $e \defeq s_1^* s_1 s_2^* s_2$. This projection $e$ is of the form $e = e_X$ for some $X \in \cJ$. We now claim that $\lambda(s_1 e) = \lambda(s_2 e)$. First of all, we have $\lambda(e) = E_X$. For $x \in X$, $\lambda(s_1 e) \ve_x = \ve_{gx}$ as $e$ is dominated by the support projection $s_1^* s_1$ of $s_1$. Similarly, $\lambda(s_2 e) \ve_x = \ve_{gx}$ for all $x \in X$. Since we clearly have $\lambda(s_1 e) \ve_y = \lambda(s_2 e) \ve_y = 0$ for $y \notin X$, we have shown $\lambda(s_1 e) = \lambda(s_2 e)$. But $\lambda$ is injective on $S$ by Lemma~\ref{lambda-inj}, so that $s_1 e = s_2 e$. Hence by definition of the equivalence relation on $\hat{S}$, we conclude that $[s_1,\chi_1] = [s_2,\chi_2]$.

To prove surjectivity of $\Phi$, take $(\chi,g) \in \cG_N^N$. $r(\chi,g) = \chi \in N$ and $d(\chi,g) = \chi g \in N$ imply that $\chi(E_{P \cap (g \cdot P)}) = \chi(E_P) (\chi g)(E_P) = 1$. As $\ping$ is Toeplitz, there exists $s \in S$ such that $E_P \lambda_g E_P = \lambda(s)$. Thus $(\chi g)(\lambda(s^* s)) = \chi(\tau_g(E_P \lambda_{g^{-1}} E_P \lambda_g E_P)) = \chi(\tau_g(E_{(g^{-1} \cdot P) \cap P})) = \chi(E_{P \cap (g \cdot P)}) = 1$. Thus $(s,(\chi g)\vert_{D_r})$ lies in $\hat{S}$. Since $g(s)=g$, we obtain $\Phi([s,(\chi g)\vert_{D_r}]) = (c^* ((\chi g)\vert_{D_r}) g^{-1},g) = (\chi,g)$.

Let us now prove that $\Phi$ is compatible with the groupoid structures. $[s_1,\chi_1]$ and $[s_2,\chi_2]$ are composable if and only if
\bgln
\label{compcond}
  s_2.\chi_2 = \chi_1 \nonumber &\LRarr& c^*\chi_1 = c^*(s_2.\chi_2) \overset{\eqref{cschi}}{=} (c^*\chi_2) g(s_2)^{-1} \nonumber \\
  &\LRarr& ((c^*\chi_1)g(s_1)^{-1})g(s_1) = (c^*\chi_2) g(s_2)^{-1}.
\egln
But this last equation is precisely the condition for composability of $\Phi([s_1,\chi_1]) = ((c^*\chi_1)g(s_1)^{-1},g(s_1))$ and $\Phi([s_2,\chi_2]) = ((c^*\chi_2)g(s_2)^{-1},g(s_2))$. If \eqref{compcond} is satisfied, then
\bglnoz
  && \Phi([s_1,\chi_1][s_2,\chi_2]) = \Phi([s_1s_2,\chi_2]) \\
  &=& ((c^*\chi_2)g(s_1s_2)^{-1},g(s_1s_2)) = ((c^*\chi_2)g(s_2)^{-1}g(s_1)^{-1},g(s_1)g(s_2)) \\
  &\overset{\eqref{compcond}}{=}& ((c^*\chi_1)g(s_1)^{-1},g(s_1)g(s_2)) = ((c^*\chi_1)g(s_1)^{-1},g(s_1)) ((c^*\chi_2)g(s_2)^{-1},g(s_2)) \\
  &=& \Phi([s_1,\chi_1])\Phi([s_2,\chi_2]).
\eglnoz
Moreover,
\bglnoz
  && \Phi([s,\chi]^{-1}) = \Phi([s^*,s.\chi]) = (c^*(s.\chi) g(s^*)^{-1}, g(s^*)) \\
  &\overset{\eqref{cschi}}{=}& ((c^*\chi)g(s)^{-1} g(s), g(s)^{-1}) = (c^*\chi,g(s)^{-1}) \\
  &=& ((c^*\chi)g(s)^{-1},g(s))^{-1} = (\Phi([s,\chi]))^{-1}.
\eglnoz

Finally, $\Phi$ is continuous by definition of the quotient topology as $\hat{S} \to \cG_N^N$, $(s,\chi) \ma ((c^*\chi)g(s)^{-1},g(s))$ is continuous. In addition, $\Phi$ is open as well. Namely, let $\pi: \hat{S} \to \cG(S)$ be the canonical projection, and take an open subset $U$ of $\cG(S)$. Then $\pi^{-1}(U)$ is open in $\hat{S}$. As $\hat{S} = \bigcup_{s \in S} \gekl{s} \times \Omega$, we must have that $\pi^{-1}(U) \cap (\gekl{s} \times \Omega)$ is open for every $s \in S$. In other words, for every $s$ in $S$ there exists an open subset $U_s$ of $\Omega$ such that $\pi^{-1}(U) = \bigcup_{s \in S} \gekl{s} \times U_s$. Hence $\Phi(U)$ = $\bigcup_{s \in S} (c^*(U_s) g(s)^{-1}) \times \gekl{g(s)}$ is open in $\cG$ as $c^*$ is open by Corollary~\ref{c-open}.
\eproof

\bremark
In particular, Proposition~\ref{id-groupoids} shows that $\cG(S)$ is Hausdorff.
\eremark

\subsection{Groupoid crossed products}

We follow \cite{Qui-Sie} and \cite{Kho-Ska2} and describe $(A \otimes D) \rtimes_{\beta} S$ as a groupoid crossed product by $\cG(S)$.

First of all, we think of $S$ as a subsemigroup of $\cG(S)^{\text{op}}$, the inverse semigroup of open $\cG(S)$-sets, via the embedding
\bgloz
  S \to \cG(S)^{\text{op}}, \: s \ma O_s \defeq \pi(\gekl{s} \times \menge{\chi \in \Omega}{\chi(\lambda(s^*s))=1})
\egloz
where $\pi$ is the canonical projection $\pi: \hat{S} \to \cG(S)$. This embedding is explained in \cite{Kho-Ska2}, directly after Theorem~6.5.

Let us now define a groupoid dynamical system $(A \times \Omega,\cG(S),\alpha(S))$ in the sense of \cite{Muh-Wil}, \S~4.1. We let $A \times \Omega$ be the trivial C*-bundle over $\Omega$ with constant fibres $A$. Consider for every $[s,\chi] \in \cG(S)$ the automorphism
\bgloz
  \alpha(S)_{[s,\chi]}: A \times \gekl{\chi} \to A \times \gekl{s.\chi}, \: (a,\chi) \ma (\alpha_{g(s)}(a),s.\chi).
\egloz
It is straightforward to check that this family $(\alpha(S)_{[s,\chi]})_{[s,\chi] \in \cG(S)}$ gives rise to the desired groupoid dynamical system in the sense of \cite{Muh-Wil}, \S~4.1. Moreover, it is also easy to see that the dynamical systems $(A \otimes D,S,\beta)$ and $(A \times \Omega,\cG(S),\alpha(S))$ correspond to one another in the sense of \cite{Qui-Sie}, Theorem~5.3. In such a situation, we may apply Theorem~7.2 of \cite{Qui-Sie} and deduce
\bprop
\label{cd2}
The map
$
  (a \otimes sds^*) \delta_s 
  \ma 
  \eckl{[t,\psi]
  \ma 
  \bfa
    \psi(\lambda(d))a \falls [t,\psi] \in O_s \\
    0 \sonst
  \efa
  }
$
extends to an isomorphism $(A \otimes D) \rtimes_{\beta} S \cong (A \times \Omega) \rtimes_{\alpha(S)} \cG(S)$.
\eprop

To proceed, we describe the full and reduced crossed products of $(A \otimes D_P^G,G,\alpha \otimes \tau)$ as groupoid crossed products. We just have to follow Example~4.8 in \cite{Muh-Wil} and \S~6 of \cite{Sims-Wil}.

The action of the transformation groupoid $\cG = \Omega_P^G \rtimes G$ on the trivial C*-bundle $A \times \Omega_P^G$ over $\Omega_P^G$ is given by the automorphisms $(\alpha \times \Omega_P^G)_{(\chi,g)}: A \times \gekl{\chi g} \to A \times \gekl{\chi}$, $(a,\chi g) \ma (\alpha_g(a),\chi)$. Identifying $A \otimes D_P^G$ with $C_0(\Omega_P^G,A)$ in the canonical way, we obtain from Example~4.8 in \cite{Muh-Wil} and \S~6 of \cite{Sims-Wil}:
\bprop
\label{cd5}
The map $C_c(G,A \otimes D_P^G) \ni f \ma \eckl{(\chi,g) \ma f(g)(\chi)} \in (A \times \Omega_P^G) \rtimes_{\alpha \times \Omega_P^G} \cG$ extends to an isomorphism $(A \otimes D_P^G) \rtimes_{\alpha \otimes \tau} G \cong (A \times \Omega_P^G) \rtimes_{\alpha \times \Omega_P^G} \cG$.

Similarly, the map $C_c(G,A \otimes D_P^G) \ni f \ma \eckl{(\chi,g) \ma f(g)(\chi)} \in (A \times \Omega_P^G) \rtimes_{\alpha \times \Omega_P^G,r} \cG$ extends to an isomorphism $(A \otimes D_P^G) \rtimes_{\alpha \otimes \tau,r} G \cong (A \times \Omega_P^G) \rtimes_{\alpha \times \Omega_P^G,r} \cG$.
\eprop

Let us now restrict the $\cG$-action $\alpha \times \Omega_P^G$ to $\cG_N^N$. We obtain an action $\alpha \times N$ of $\cG_N^N$ on the sub-C*-bundle $A \times N$ (i.e. just the trivial C*-bundle over $N$ with constant fibres $A$). The following observation links the two groupoid dynamical systems we are considering:
\blemma
The dynamical systems $(A \times \Omega,\cG(S),\alpha(S))$ and $(A \times N,\cG_N^N,\alpha \times N)$ are isomorphic. More precisely, the identifications $\id \times c^*: A \times \Omega \cong A \times N$ and $\Phi: \cG(S) \to \cG_N^N$, $[s,\chi] \ma ((c^*\chi)g(s)^{-1},g(s))$ transport the action $\alpha(S)$ to $\alpha \times N$, in the sense that for every $[s,\chi] \in \cG(S)$ and $(a,\chi) \in A \times \Omega$, we have $(\alpha \times N)_{\Phi([s,\chi])}((\id \times c^*)(a,\chi)) = (\id \times c^*)(\alpha(S)_{[s,\chi]}(a,\chi))$.
\elemma
\bproof
We just have to compute that
\bglnoz
  && (\alpha \times N)_{\Phi([s,\chi])}((\id \times c^*)(a,\chi)) = (\alpha \times N)_{((c^*\chi)g(s)^{-1},g(s))}(a,c^*\chi) \\
  &=& (\alpha_{g(s)}(a),(c^*\chi)g(s)^{-1}) \overset{\eqref{cschi}}{=} (\alpha_{g(s)}(a),c^*(s.\chi)) \\
  &=& (\id \times c^*)(\alpha_{g(s)}(a),s.\chi) = (\id \times c^*)(\alpha(S)_{[s,\chi]}(a,\chi)).
\eglnoz
\eproof
\bcor
\label{cd3}
The map $C_c(\cG(S),A) \ni f \ma f \circ \Phi^{-1} \in C_c(\cG_N^N,A)$ extends to isomorphisms
\bglnoz
  && (A \times \Omega) \rtimes_{\alpha(S)} \cG(S) \overset{\cong}{\lori} (A \times N) \rtimes_{\alpha \times N} \cG_N^N, \\
  && (A \times \Omega) \rtimes_{\alpha(S),r} \cG(S) \overset{\cong}{\lori} (A \times N) \rtimes_{\alpha \times N,r} \cG_N^N.
\eglnoz
\ecor

We now want to see that the $(\cG,\cG_N^N)$-equivalence $\cG_N$ of Corollary~\ref{GGNN} gives rise to an equivalence between $(A \times \Omega_P^G,\cG,\alpha \times \Omega_P^G)$ and $(A \times N,\cG_N^N,\alpha \times N)$.

\blemma
Equip the trivial Banach-bundle $A \times \cG_N$ with the fibrewise imprimitivity bimodule structure given by the inner products
\bglnoz
  && {}_{A \times \gekl{r(\gamma)}} \spkl{(a_1,\gamma),(a_2,\gamma)} = (a_1 a_2^*,r(\gamma)) \in A \times \Omega_P^G, \\
  && \spkl{(a_1,(\chi,g)),(a_2,(\chi,g))}_{A \times \gekl{d(\gamma)}} = (\alpha_g^{-1}(a_1^* a_2),\chi g) \in A \times N
\eglnoz
and the left and right actions
\bglnoz
  && (a_l,r(\gamma)) \cdot (a,\gamma) = (a_l a,\gamma) \text{ for } (a_l,r(\gamma)) \in A \times \Omega_P^G, \\
  && (a,(\chi,g)) \cdot (a_r,\chi g) = (a \alpha_g(a_r),(\chi,g)) \text{ for } (a_r,\chi g) \in A \times N.
\eglnoz
Moreover, let $\cG$ act from the left on $A \times \cG_N$ by $(\chi_l,g_l) \cdot (a,\gamma) = (\alpha_{g_l}(a),(\chi_l,g_l)\gamma)$ and let $\cG_N^N$ act from the right on $A \times \cG_N$ by $(a,\gamma) \cdot \gamma_r = (a,\gamma \gamma_r)$.

Then in this way, $A \times \cG_N$ becomes an equivalence between $(A \times \Omega_P^G,\cG,\alpha \times \Omega_P^G)$ and $(A \times N,\cG_N^N,\alpha \times N)$ in the sense of Definition~5.1 in \cite{Muh-Wil}.
\elemma
\bproof
Just verify by straightforward computations that the axioms for an equivalence in Definition~5.1 from \cite{Muh-Wil} are satisfied.
\eproof

The reason why this is interesting for us is the following consequence of Theorem~5.5 in \cite{Muh-Wil} and Corollary~19 from \cite{Sims-Wil}:
\blemma
The canonical inclusion $C_c(\cG_N^N,A) \into C_c(\cG,A)$ extends to (isometric!) embeddings
\bglnoz
  && (A \times N) \rtimes_{\alpha \times N} \cG_N^N \into (A \times \Omega_P^G) \rtimes_{\alpha \times \Omega_P^G} \cG, \\
  && (A \times N) \rtimes_{\alpha \times N,r} \cG_N^N \into (A \times \Omega_P^G) \rtimes_{\alpha \times \Omega_P^G,r} \cG.
\eglnoz
\elemma
\bproof
As $N$ is clopen, $\cG_N^N$ is a clopen subset of $\cG$, so that we really have $C_c(\cG_N^N,A) \subseteq C_c(\cG,A)$. As $\cG_N$ is clopen as well, we actually have $C_c(\cG_N^N,A) \subseteq C_c(\cG_N,A) \subseteq C_c(\cG,A)$.

Let us first treat the full crossed products. Take a function $f \in C_c(\cG_N^N,A)$. All we have to show is that
\bgl
\label{toshow}
  \norm{f}_{(A \times N) \rtimes_{\alpha \times N} \cG_N^N} = \norm{f}_{(A \times \Omega_P^G) \rtimes_{\alpha \times \Omega_P^G} \cG}.
\egl
We denote by $*$ the convolution product in $(A \times \Omega_P^G) \rtimes_{\alpha \times \Omega_P^G} \cG$, and observe that its restriction to $C_c(\cG_N^N,A)$ coincides with the convolution product coming from $(A \times N) \rtimes_{\alpha \times N} \cG_N^N$. We certainly have
\bgln
\label{crucial1}
  && \norm{f}_{(A \times N) \rtimes_{\alpha \times N} \cG_N^N}^2 = \norm{f^* * f}_{(A \times N) \rtimes_{\alpha \times N} \cG_N^N}, \nonumber \\
  && \norm{f}_{(A \times \Omega_P^G) \rtimes_{\alpha \times \Omega_P^G} \cG}^2 = \norm{f * f^*}_{(A \times \Omega_P^G) \rtimes_{\alpha \times \Omega_P^G} \cG}.
\egln
But now comes the crucial observation, namely that
\bgl
\label{crucial}
  f_1^* * f_2 = \spkl{\spkl{f_1,f_2}}_{(A \times N) \rtimes_{\alpha \times N} \cG_N^N} 
  \text{ and }
  f_1 * f_2^* = {}_{(A \times \Omega_P^G) \rtimes_{\alpha \times \Omega_P^G} \cG} \spkl{\spkl{f_1,f_2}}
\egl
for all $f_1$, $f_2$ in $C_c(\cG_N,A)$. Here $\spkl{\spkl{\cdot,\cdot}}$ are the inner products defined in Theorem~5.5 in \cite{Muh-Wil}. The verification of \eqref{crucial} is a straightforward computation. In particular, we have for our function $f$
\bgl
\label{crucial2}
  f^* * f = \spkl{\spkl{f,f}}_{(A \times N) \rtimes_{\alpha \times N} \cG_N^N} 
  \text{ and }
  f * f^* = {}_{(A \times \Omega_P^G) \rtimes_{\alpha \times \Omega_P^G} \cG} \spkl{\spkl{f,f}}.
\egl
By Theorem~5.5 in \cite{Muh-Wil}, $C_c(\cG_N,A)$ is a $(A \times \Omega_P^G) \rtimes_{\alpha \times \Omega_P^G} \cG$-$(A \times N) \rtimes_{\alpha \times N} \cG_N^N$ pre-imprimitivity bimodule. Therefore, we conclude that
\bgl
\label{crucial3}
  \norm{\spkl{\spkl{f,f}}_{(A \times N) \rtimes_{\alpha \times N} \cG_N^N}} 
  = \norm{{}_{(A \times \Omega_P^G) \rtimes_{\alpha \times \Omega_P^G} \cG} \spkl{\spkl{f,f}}},
\egl
where we take the norm in $(A \times N) \rtimes_{\alpha \times N} \cG_N^N$ on the left hand side and the norm in $(A \times \Omega_P^G) \rtimes_{\alpha \times \Omega_P^G} \cG$ on the right hand side. Inserting \eqref{crucial2} and \eqref{crucial1} into \eqref{crucial3}, we obtain \eqref{toshow}, as desired.

To treat reduced crossed products, just use Theorem~14 or rather Corollary~19 of \cite{Sims-Wil} instead of Theorem~5.5 in \cite{Muh-Wil}.
\eproof

\bcor
\label{cd4}
The inclusion $C_c(\cG_N^N,A) \into C_c(\cG,A)$ extends to isomorphisms
\bglnoz
  && (A \times N) \rtimes_{\alpha \times N} \cG_N^N \cong \1z_N \rukl{(A \times \Omega_P^G) \rtimes_{\alpha \times \Omega_P^G} \cG} \1z_N, \\
  && (A \times N) \rtimes_{\alpha \times N,r} \cG_N^N \cong \1z_N \rukl{(A \times \Omega_P^G) \rtimes_{\alpha \times \Omega_P^G,r} \cG} \1z_N.
\eglnoz
Here $\1z_N$ is the characteristic function of $N \subseteq \cG$, viewed in a canonical way as a multiplier of $(A \times \Omega_P^G) \rtimes_{\alpha \times \Omega_P^G} \cG$ and $(A \times \Omega_P^G) \rtimes_{\alpha \times \Omega_P^G,r} \cG$, respectively. Moreover, $\1z_N \rukl{(A \times \Omega_P^G) \rtimes_{\alpha \times \Omega_P^G} \cG} \1z_N$ and $\1z_N \rukl{(A \times \Omega_P^G) \rtimes_{\alpha \times \Omega_P^G,r} \cG} \1z_N$ are full corners in the corresponding full and reduced crossed products.
\ecor
\bproof
It is easy to see that $C_c(\cG_N^N,A) = \1z_N * C_c(\cG,A) * \1z_N$. Thus the first part of the corollary follows from the previous lemma. We also have $C_c(\cG_N,A) = C_c(\cG,A) * \1z_N$. Using this observation and also \eqref{crucial} from the proof of the previous lemma, the second part of our assertion follows from \cite{Muh-Wil}, Theorem 5.5 in the case of full crossed products and from \cite{Sims-Wil}, Corollary~19 in the reduced case.
\eproof

Let us summarize what we have obtained so far:
\btheo
\label{thm1}
Let $P$ be a subsemigroup of a group $G$. Assume that $\cJ$ is independent (see Definition~\ref{ind}) and that $\ping$ satisfies the Toeplitz condition from Definition~\ref{Toeplitz}. Then the following diagram commutes:
\bgl
\label{cd}
  \begin{CD}
  A \rta_{\alpha,s} P @> \lambda_{(A,P,\alpha)} >> A \rta_{\alpha,r} P \\
  @VV \cong V @V \cong VV \\
  (A \otimes D) \rtimes_{\beta} S @>>> (A \otimes D) \rtimes_{\beta,r} S \\
  @VV \cong V @V \cong VV \\
  (A \times \Omega) \rtimes_{\alpha(S)} \cG(S) @>>> (A \times \Omega) \rtimes_{\alpha(S),r} \cG(S) \\
  @VV \cong V @V \cong VV \\
  (A \times N) \rtimes_{\alpha \times N} \cG_N^N @>>> (A \times N) \rtimes_{\alpha \times N,r} \cG_N^N \\
  @VV \cong V @V \cong VV \\
  \1z_N \rukl{(A \times \Omega_P^G) \rtimes_{\alpha \times \Omega_P^G} \cG} \1z_N
  @>>> 
  \1z_N \rukl{(A \times \Omega_P^G) \rtimes_{\alpha \times \Omega_P^G,r} \cG} \1z_N \\
  @VV \cong V @V \cong VV \\
  (1 \otimes E_P) \rukl{(A \otimes D_P^G) \rtimes_{\alpha \otimes \tau} G} (1 \otimes E_P)
  @>>> 
  (1 \otimes E_P) \rukl{(A \otimes D_P^G) \rtimes_{\alpha \otimes \tau,r} G} (1 \otimes E_P)
  \end{CD}
\egl

Moreover, $\1z_N$ and $1 \otimes E_P$ give rise to full corners in the full and reduced crossed products associated with $(A \times \Omega_P^G,\cG,\alpha \times \Omega_P^G)$ and $(A \otimes D_P^G,G,\alpha \otimes \tau)$.

And finally, the square at the bottom of diagram~\eqref{cd} is obtained by restricting the commutative diagram
\bgl
\label{cd-group}
  \begin{CD}
  (A \times \Omega_P^G) \rtimes_{\alpha \times \Omega_P^G} \cG @>>> (A \times \Omega_P^G) \rtimes_{\alpha \times \Omega_P^G,r} \cG \\
  @VV \cong V @V \cong VV \\
  (A \otimes D_P^G) \rtimes_{\alpha \otimes \tau} G @>>> (A \otimes D_P^G) \rtimes_{\alpha \otimes \tau,r} G
  \end{CD}
\egl
In all these diagrams, the horizontal arrows are given by the canonical projections (the regular representations), and the vertical maps are the isomorphisms we have explicitly constructed before.
\etheo
\bproof
We just have to collect what we have proven. The first commuting square from the top is given by Proposition~\ref{cd1} and Lemma~\ref{reduced-invsemigp-cropro}. Proposition~\ref{cd2} tells us that on the left hand side, the second vertical arrow from the top is an isomorphism. The third square and its commutativity is provided by Corollary~\ref{cd3}. The fourth square and its commutativity is given by Corollary~\ref{cd4}. And Proposition~\ref{cd5} gives the square at the bottom of \eqref{cd} and that it is the restriction of the commutative diagram~\eqref{cd-group}. Using Corollary~\ref{full-red}, we can fill in the vertical arrow on the right hand side of the second square so that it commutes as well. That $\1z_N$ gives rise to full corners is shown in Corollary~\ref{cd4}, and it corresponds to $1 \otimes E_P$ under the isomorphism from Proposition~\ref{cd5}. This completes the proof.
\eproof

\section{Nuclearity}
\label{nuc}

Using our results from the previous section, we obtain equivalent characterizations for nuclearity of semigroup C*-algebras.

\btheo
\label{thm2}
Let $P$ be a subsemigroup of a group $G$. Assume that $\cJ$ is independent (see Definition~\ref{ind}) and that $\ping$ satisfies the Toeplitz condition from Definition~\ref{Toeplitz}. Then the following are equivalent:
\begin{enumerate}
\item[(i)] $C^*_s(P)$ is nuclear.
\item[(ii)] $C^*_r(P)$ is nuclear.
\item[(iii)] Whenever given a $G$-action $\alpha$ on a C*-algebra $A$, the canonical homomorphism $\lambda_{(A,P,\alpha)}: A \rta_{\alpha,s} P \to A \rta_{\alpha,r} P$ is an isomorphism.
\item[(iv)] The groupoid $\cG_N^N$ is amenable.
\item[(v)] The groupoid $\cG$ is amenable.
\end{enumerate}

Here $\cG$ is the transformation groupoid $\Omega_P^G \rtimes G$ from \S~\ref{groupoids}, and the groupoid $\cG_N^N$ is the restriction of $\cG$ also introduced in \S~\ref{groupoids}.
\etheo
Of course, amenability of $\cG$ just means that $G$ acts amenably on $\Omega_P^G$.

\bproof
We prove \an{(ii) $\LRarr$ (iv) $\LRarr$ (v)} and \an{(iv) $\Rarr$ (iii) $\Rarr$ (i) $\Rarr$ (ii)}.

To see \an{(ii) $\LRarr$ (iv)}, plug in $A = \Cz$ in diagram~\eqref{cd}. Moreover, using that $(\Cz \times N) \rtimes_{\tr \times N,r} \cG_N^N$ is canonically isomorphic to the reduced groupoid C*-algebra $C^*_r(\cG_N^N)$ by Example~10 in \cite{Sims-Wil}, we see that $C^*_r(P) \cong C^*_r(\cG_N^N)$. Since $\cG_N^N$ is an $r$-discrete (also called \'{e}tale) groupoid, it is known that $C^*_r(\cG^N_N)$ is nuclear if and only if $\cG_N^N$ is amenable (see for instance \cite{Br-Oz}, Chapter~5, Theorem~6.18).

For \an{(iv) $\LRarr$ (v)}, recall that we have proven that $\cG$ is equivalent to $\cG_N^N$ in Corollary~\ref{GGNN}. As amenability is invariant under equivalences of groupoids by Theorem~2.2.17 in \cite{An-Ren}, we have proven \an{(iv) $\LRarr$ (v)}.

To see \an{(iv) $\Rarr$ (iii)}, note that amenability of $\cG_N^N$ implies that the fourth (counted from the top) horizontal map in diagram~\eqref{cd} is an isomorphism by \cite{An-Ren}, Proposition~6.1.10. By commutativity of \eqref{cd}, we deduce that $\lambda_{(A,P,\alpha)}$ must be an isomorphism.

For \an{(iii) $\Rarr$ (i)}, first apply (iii) to $A = \Cz$ to deduce that $\lambda: C^*_s(P) \to C^*_r(P)$ is an isomorphism. Now use an argument as in the proof of Theorem~4.3, 3 $\Rarr$ 4 in \cite{Br-Oz}, Chapter~4: By the definition of semigroup crossed products by automorphisms, it is easily seen that for trivial actions, we can canonically identify $A \rta_{\tr,s} P$ with $A \otimes_{\max} C^*_s(P)$ and $A \rta_{\tr,r} P$ with $A \otimes_{\min} C^*_r(P) \overset{\id \otimes_{\min} \lambda^{-1}}{\cong} A \otimes_{\min} C^*_s(P)$ such that the diagram
\bgloz
  \begin{CD}
  A \rta_{\tr,s} P @> \lambda_{(A,P,\alpha)} >> A \rta_{\tr,r} P \\
  @VV \cong V @V \cong VV \\
  A \otimes_{\max} C^*_s(P) @>>> A \otimes_{\min} C^*_s(P)
  \end{CD}
\egloz
commutes. The horizontal map at the bottom is the canonical homomorphism, and it must be an isomorphism since $\lambda_{(A,P,\alpha)}$ is one by (iii). This means that $C^*_s(P)$ is nuclear.

Finally, to go from (i) to (ii), just observe that $C^*_r(P)$ is a quotient of $C^*_s(P)$ and apply \cite{Bla}, Corollary~IV.3.1.13.
\eproof

\bremark
\label{nuc-am}
In particular, we see that in the situation of Theorem~\ref{thm2}, nuclearity of $C^*_s(P)$ (or $C^*_r(P)$) implies that the left regular representation $\lambda: C^*_s(P) \to C^*_r(P)$ is faithful.
\eremark

\bremark
\label{unital-suff}
It certainly suffices to consider unital $A$ in (iii) of Theorem~\ref{thm2}.
\eremark

For later purposes, we derive the following consequence:
\bcor
\label{Theta}
Let $P$ be a subsemigroup of a group $G$. Assume that $\cJ$ is independent and that $\ping$ satisfies the Toeplitz condition. If $C^*_r(P)$ is nuclear, then there exists a net of completely positive contractions $\Theta_i: C^*_r(P) \to C^*_r(P)$ such that
\begin{enumerate}
\item[1.] $\lim_i \Theta_i(x) = x$ for all $x \in C^*_r(P)$, 
\item[2.] for every $i$ there is a finitely supported function $d_i: G \to D_r$, $g \ma d_i(g)$ such that $\Theta_i(V) = d_i(g_r(V))V$ for all $0 \neq V \in S_r$. 
\end{enumerate}
($S_r$ and $g_r$ were introduced in \S~\ref{invsemi} and Lemma~\ref{g}.)
\ecor
\bproof
As $C^*_r(P)$ is nuclear, $\cG = \Omega_P^G \rtimes G$ is amenable by the previous theorem. Thus combining Theorem~6.18 and Proposition~6.16 from Chapter~5 in \cite{Br-Oz}, we obtain a net $h_i \in C_c(\cG)$ such that $h_i \lori_i 1$ uniformly on compact subsets and $C_c(\cG) \ni f \ma h_i \cdot f \in C_c(\cG)$ extends to a completely positive contraction on $C^*_r(\cG)$. Under the canonical identification $D_P^G \rtimes_{\tau,r} G \cong C^*_r(\cG)$, we obtain a net of completely positive contractions $m_i$ such that $m_i(x) \lori_i x$ for all $x \in D_P^G \rtimes_{\tau,r} G$ and for every $i$, there exists a finitely supported function $\ti{h}_i: G \to D_P^G$ with $m_i(d U_g) = \ti{h}_i(g) d U_g$ for all $d \in D_P^G$ and $g \in G$. To be precise, $\ti{h}_i(g)$ is given by $\chi \ma h_i(\chi,g)$. Now, let $\Theta_i$ be the composition $C^*_r(P) \cong E_P (D_P^G \rtimes_{\tau,r} G) E_P \subseteq D_P^G \rtimes_{\tau,r} G \overset{m_i}{\lori} D_P^G \rtimes_{\tau,r} G$. Here we have used Lemma~\ref{AP=APG}. We have
\bglnoz
  && \Theta_i(V_{p_1}^* V_{q_1} \dotsm V_{p_n}^* V_{q_n}) = m_i(E_{p_1^{-1} q_1 \dotsm p_n^{-1} q_n P} U_{p_1^{-1} q_1 \dotsm p_n^{-1} q_n}) \\
  &=& (E_P \ti{h}_i(p_1^{-1} q_1 \dotsm p_n^{-1} q_n) E_P) E_{p_1^{-1} q_1 \dotsm p_n^{-1} q_n P} U_{p_1^{-1} q_1 \dotsm p_n^{-1} q_n}.
\eglnoz
Set $d_i(g) = E_P \ti{h}_i(g) E_P$. Then $d_i$ lies in $D_r$ by the Toeplitz condition. Moreover, we see that $\Theta_i$ has image in $E_P (D_P^G \rtimes_{\tau,r} G) E_P$, so that identifying this corner back again with $C^*_r(P)$, we obtain the desired net of completely positive contractions.
\eproof
This observation will be used in the next section when we study induced ideals of semigroup C*-algebras. Now let us show that the existence of such completely positive contractions $\Theta_i$ on $C^*_s(P)$ implies nuclearity of $C^*_s(P)$. First, we set $D_g \defeq \lspan(\menge{s \in S}{g(s)=g}) \subseteq C^*_s(P)$, and for a map $\Theta$ on $C^*_s(P)$, we let the d-support of $\Theta$ be $\dsupp(\Theta) = \menge{g \in G}{\Theta \vert_{D_g} \neq 0}$. By Theorem~\ref{thm2}, we know that under the assumptions of Corollary~\ref{Theta}, $\lambda: C^*_s(P) \to C^*_r(P)$ is an isomorphism. Thus Corollary~\ref{Theta} gives us completely positive contractions $\Theta_i$ on $C^*_s(P)$ ($\cong C^*_r(P)$) such that $\lim_i \Theta_i(x) = x$ for all $x \in C^*_s(P)$ and $\abs{\dsupp(\Theta_i)} < \infty$ for all $i$. The following result shows that the existence of such $\Theta_i$ already implies nuclearity of $C^*_s(P)$:
\bprop
\label{Theta->nuc}
Let $P$ be a subsemigroup of a group $G$. Assume that $\cJ$ is independent. Moreover, assume that there exists a net of completely positive contractions $\Theta_i: C^*_s(P) \to C^*_s(P)$ such that
\bgln
\label{w*}
  && \lim_i \Theta_i(x) = x \fa x \in C^*_s(P), \\
\label{finite-d-supp}
  && \abs{\dsupp(\Theta_i)} < \infty \fa i.
\egln
Then for every C*-algebra $A$, $\lambda_{(A,P,\tr)}: A \rta_{\tr,s} P \to A \rta_{\tr,r} P$ is an isomorphism.

In particular, $\lambda: C^*_s(P) \to C^*_r(P)$ is an isomorphism and $C^*_s(P)$ is nuclear.
\eprop
Note that we do not assume that $\ping$ is Toeplitz.
\bproof
The proof is just the same as the one for \an{5) $\Rarr$ 6)} in \cite{Li2}, \S~4, but for arbitrary coefficients. For the sake of completeness, we write out the proof.


Let $\cE_s^A$ be the composite $A \rta_{\tr,s} P \overset{\lambda_{(A,P,\tr)}}{\lori} A \rtaar P \cong A \otimes_{\min} C^*_r(P) \overset{\id \otimes \cE_r}{\lori} A \otimes D_r \overset{\id \otimes (\lambda \vert_D)^{-1}}{\lori} A \otimes D \to A \rta_{\tr,s} P$. Here $\cE_r$ is the conditional expectation $C^*_r(P) \to D_r$ from \cite{Li2}, \S~3.2. Moreover, we used that $\lambda \vert_D: D \to D_r$ is an isomorphism as $\cJ$ is independent. The last homomorphism is given by $A \otimes D \to A \rta_{\tr,s} P$, $a \otimes d \ma \iota(a) \overline{d}$. Now set $D_g^A \defeq \lspan(\menge{\iota(a) \overline{x}}{a \in A, x \in D_g}) \subseteq A \rta_{\tr,s} P$. Obviously the algebraic sum $\sum_{g \in G} D_g^A$ is dense in $A \rta_{\tr,s} P$. Moreover, we have by construction that $\cE_s^A \vert_{D_e^A} = \id_{D_e^A}$ and $\cE_s^A \vert_{D_g^A} = 0$ if $g \neq e$.

Given a positive functional $\phi$ on $A \rtaas P$, set $\dsupp(\phi) = \menge{g \in G}{\phi \vert_{D_g^A} \neq 0}$. If $\abs{\dsupp(\phi)} < \infty$, then we have for all $x \in A \rta_{\tr,s} P$:
\bgl
\label{phi-E}
  \abs{\phi(x)}^2 \leq \abs{\dsupp(\phi)} \norm{\phi} \phi(\cE_s^A(x^*x)).
\egl
To prove \eqref{phi-E}, let $\dsupp(\phi) = \gekl{g_1, \dotsc, g_n}$. As $\sum_{g \in G} D_g^A$ is dense in $A \rta_{\tr,s} P$, it suffices to prove \eqref{phi-E} for $x \in \sum_{g \in G} D_g^A$. Take such an element $x$ and a finite subset $F \subseteq G$ such that $\dsupp(\phi) \subseteq F$ and $x = \sum_{g \in F} x_g$ with $x_g \in D_g^A$. Then the same computation as in the proof of Lemma~4.8 in \cite{Li2} yields
\bglnoz
  \abs{\phi(x)}^2 &=& \abs{\sum_{i=1}^n \phi(x_{g_i})}^2 \leq n \sum_{i=1}^n \abs{\phi(x_{g_i})}^2 \leq n \norm{\phi} \sum_{i=1}^n \phi(x_{g_i}^* x_{g_i}) \\
  &\leq& n \norm{\phi} \phi(\sum_{g \in F} x_g^* x_g) = n \norm{\phi} \phi(\cE_s^A(\sum_{g,h \in F} x_g^* x_h)) \\
  &=& \abs{\dsupp(\phi)} \norm{\phi} \phi(\cE_s^A(x^*x)).
\eglnoz
This proves \eqref{phi-E}.

Now take $x \in \ker(\lambda_{(A,P,\tr)})$, $x \geq 0$, and a positive functional $\phi$ on $A \rta_{\tr,s} P$. Let $\phi_i$ be the composition $A \rta_{\tr,s} P \cong A \otimes_{\max} C^*_s(P) \overset{\id \otimes \Theta_i}{\lori} A \otimes_{\max} C^*_s(P) \cong A \rta_{\tr,s} P \overset{\phi}{\to} \Cz$. These positive functionals $\phi_i$ satisfy $\lim_i \phi_i(x) = \phi(x)$ and $\abs{\dsupp(\phi_i)} < \infty$. As $\lambda_{(A,P,\tr)}(x) = 0$, we must have $\cE_s^A(x^*x) = 0$ by construction of $\cE_s^A$. Thus by \eqref{phi-E}, we conclude that $\phi_i(x) = 0$ for all $i$. Therefore $\phi(x) = \lim_i \phi_i(x) = 0$. As $\phi$ was arbitrary, we conclude that $x = 0$. Hence $\lambda_{(A,P,\tr)}$ is faithful, and we have proven the first part of our proposition. To see that $\lambda: C^*_s(P) \to C^*_r(P)$ is an isomorphism, just set $A = \Cz$. And finally, to see that $C^*_s(P)$ is nuclear, just proceed as in the proof of Theorem~\ref{thm2}, \an{(iii) $\Rarr$ (i)}.
\eproof
Under the (rather strong) assumption of left amenability, such $\Theta_i$ always exist:
\blemma
If $P$ is cancellative and left amenable, then there exists a net $\Theta_i$ as in Proposition~\ref{Theta->nuc} satisfying \eqref{w*} and \eqref{finite-d-supp}.
\elemma
\bproof
First of all, $P$ embeds into a group if it is cancellative and left amenable (see for instance \cite{Li2}, Corollary~4.5), so that we can form $C^*_s(P)$. As 5) in \cite{Li2}, \S~4.1 holds if $P$ is left cancellative and left amenable, we can form the states $\varphi_i: C^*_s(P) \to \Cz$, $x \ma \spkl{\lambda(x) \xi_i,\xi_i}$ with the $\xi_i$ from 5) of \S~4.1 in \cite{Li2}. Let $\Theta_i$ be the composition $C^*_s(P) \overset{\Delta}{\lori} C^*_s(P) \otimes_{\max} C^*_s(P) \overset{\varphi_i \otimes \id}{\lori} C^*_s(P)$, where $\Delta$ is given by (36) in \cite{Li2}. By construction, we have $\Theta_i(s) = \varphi_i(s) s \lori_i s$ for all $0 \neq s \in S$ by 5) in \cite{Li2}, \S~4.1. Therefore the $\Theta_i$ satisfy \eqref{w*}. As the $\xi_i$ have finite support (see \cite{Li2}, \S~4.1, 5)), it follows that $\abs{\dsupp(\Theta_i)} < \infty$ for all $i$ (compare also \cite{Li2}, \S~4.2, \an{5) $\Rarr$ 6)}).
\eproof
As a consequence, we obtain the following converse of Proposition~4.17 in \cite{Li2}:
\bcor
If $P$ is cancellative, left amenable and if $\cJ$ is independent, then $C^*_s(P)$ is nuclear.
\ecor
This result was also obtained independently in \cite{Nor} using different methods.

\section{Ideals induced from invariant spectral subsets}
\label{ideals}

In this section, we always assume that $\cJ$ is independent and that $\ping$ satisfies the Toeplitz condition. In this situation, we have seen that the full and reduced semigroup C*-algebras of $P$ can be described up to Morita equivalence as full or reduced crossed products by $G$. So in principle, this reduces questions about the ideal structure of semigroup C*-algebras to corresponding questions about certain crossed products by $G$. However, in our concrete situation, there are certain induced ideals which play a distinguished role. We first of all show that nuclearity allows us to describe induced ideals in a satisfactory way. Moreover, building on our results from Section~\ref{var-des}, we describe induced ideals (and their quotients) as crossed products by $G$ up to Morita equivalence. As these induced ideals correspond to (closed) invariant subsets of the spectrum $\Omega$ of the diagonal sub-C*-algebra $D_r$ (or $D$), we take a closer look at this spectrum. Using our observations from \S~\ref{fam-sub}, we describe it in terms of $\cJ$. Finally, we turn to the boundary of the spectrum and investigate the corresponding boundary action.

\subsection{Induced ideals}
\label{ind-ideals}

Let $I_r$ be an ideal of $D_r$, the diagonal sub-C*-algebra of $C^*_r(P)$. Restricting the canonical conditional expectation $\cL(\ell^2(P)) \to \ell^{\infty}(P)$, we obtain a conditional expectation $\cE_r: C^*_r(P) \to D_r$ (compare \cite{Li2}, \S~3.2). Following \cite{Ni1}, we define the induced ideal
\bgloz
  \Ind I_r \defeq \menge{x \in C^*_r(P)}{\Ad(V)\cE_r(x^*x) \in I_r \fa V \in S_r}.
\egloz
As A. Nica explains, the name \an{induced ideal} is justified because we could have obtained $\Ind I_r$ by an induction process as described in \cite{Ni1}, \S~6.1. 

For the purpose of inducing ideals, it suffices to consider invariant ideals of $D_r$.
\bdefin
An ideal $I_r$ of $D_r$ is called invariant if $\Ad(V)(I_r) \subseteq I_r$ for all $V \in S_r$.
\edefin
The reason why we only need to consider invariant ideals is that given an ideal $I_r$ of $D_r$, we obtain the invariant ideal $I_r^{({\rm inv})} \defeq \menge{d \in D_r}{\Ad(V)(d) \in I_r \fa V \in S_r}$. And just as in \cite{Ni1}, we have $\Ind I_r = \Ind I_r^{({\rm inv})} = \menge{x \in C^*_r(P)}{\cE_r(x^*x) \in I_r^{({\rm inv})}}$.

We observe that the induction process is an injective map from the set of invariant ideals of $D_r$ to the set of ideals of $C^*_r(P)$.
\blemma
\label{indcap}
Every invariant ideal $I_r$ of $D_r$ satisfies $(\Ind I_r) \cap D_r = \cE_r(\Ind I_r) = I_r$.
\elemma
\bproof
$(\Ind I_r) \cap D_r$ is contained in $\cE_r(\Ind I_r)$ as $\cE_r \vert_{D_r} = \id_{D_r}$. 

To see $\cE_r(\Ind I_r) \subseteq I_r$, take $x \in \Ind I_r$. Then $\cE_r(x)$ lies in $D_r$ and we have $\cE_r(x)^* \cE_r(x) \leq \cE_r(x^*x) \in I_r$. Thus $\cE_r(x)^* \cE_r(x)$ lies in $I_r$, and this implies $\cE_r(x) \in I_r$.

And finally, $I_r$ is contained in $\Ind I_r$ (as $\cE_r \vert_{I_r} = \id_{I_r}$) and in $D_r$ anyway.
\eproof

Following ideas of \cite{Ni1}, we deduce the following consequence of nuclearity:
\bprop
\label{nuc-ind-id}
If $\cJ$ is independent, if $\ping$ is Toeplitz and if $C^*_r(P)$ is nuclear, then $\Ind I_r$ coincides with the ideal $\spkl{I_r}$ of $C^*_r(P)$ generated by $I_r$.
\eprop
\bproof
It is clear that $\Ind I_r \supseteq \spkl{I_r}$ as $I_r$ is contained in $\Ind I_r$ by the previous lemma, and because $\Ind I_r$ is an ideal of $C^*_r(P)$.

To prove that $\Ind I_r \subseteq \spkl{I_r}$, first set, for $g \in G$, $(D_r)_g \defeq \clspan(\menge{V \in S_r}{g_r(V)=g})$. Moreover, let $(\Ind I_r)_c = \Ind I_r \cap (\sum_{g \in G} (D_r)_g) = \menge{x \in \sum_{g \in G} (D_r)_g}{\cE_r(x^*x) \in I_r}$. Here $\sum_{g \in G} (D_r)_g$ means the algebraic sum (without taking the closure), i.e. the set of finite sums of the form $\sum_{g \in G} x_g$ with $x_g \in (D_r)_g$.

As a first step, let us prove $(\Ind I_r)_c \subseteq \spkl{I_r}$: Take $x = \sum_g x_g \in (\Ind I_r)_c$. This means that $\cE_r(x^*x) = \sum_g x_g^* x_g$ lies in $I_r$. Hence ($I_r$ is hereditary) all the $x_g^* x_g$ lie in $I_r$. By polar decomposition (see \cite{Bla}, \S~II.3.2), we deduce that $x_g \in \spkl{I_r}$. Thus $x$ lies in $\spkl{I_r}$.

The second step is to prove $\Ind I_r \subseteq \overline{(\Ind I_r)_c}$. By Corollary~\ref{Theta}, there exists a net $\Theta_i$ of completely positive contractions $C^*_r(P) \to C^*_r(P)$ satisfying 1. and 2. from Corollary~\ref{Theta}. From 2., we deduce that for all $x \in C^*_r(P)$, we have $\cE_r(\Theta_i(x)) = d_i(e) \cE_r(x)$ as this formula obviously holds for $x \in \sum_{g \in G} (D_r)_g$ because of 2. Now take $x \in \Ind I_r$. Then $\Theta_i(x)$ lies in $\Ind I_r$ as well since $\cE_r(\Theta_i(x)^* \Theta_i(x)) \leq \cE_r(\Theta_i(x^*x)) = d_i(e) \cE_r(x^*x) \in I_r$. Moreover, as the $d_i$ in Corollary~\ref{Theta} have finite support, we have $\Theta_i(x) \in \sum_{g \in G} (D_r)_g$. Thus $\Theta_i(x)$ is in $(\Ind I_r)_c$. And finally, by 1. in Corollary~\ref{Theta}, $x = \lim_i \Theta_i(x)$ lies in $\overline{(\Ind I_r)_c}$.
\eproof

Just as in \cite{Ni1}, we obtain the following characterization of induced ideals:
\bcor
If $\cJ$ is independent, if $\ping$ is Toeplitz and if $C^*_r(P)$ is nuclear, then
\bgloz
  \menge{\Ind I_r}{I_r \triangleleft D_r} = \menge{J \triangleleft C^*_r(P)}{\cE_r(J) \subseteq J}.
\egloz
\ecor
\bproof
\an{$\subseteq$} holds by Lemma~\ref{indcap}. To prove \an{$\supseteq$}, take an ideal $J$ of $C^*_r(P)$ such that $\cE_r(J) \subseteq J$. As $J$ is an ideal of $C^*_r(P)$, $\cE_r(J)$ is an invariant ideal of $D_r$. Moreover, $J$ is contained in $\Ind \cE_r(J)$ as for $x \in J$, $x^*x$ also lies in $J$, hence $\cE_r(x^*x)$ lies in $\cE_r(J)$. Thus by the last corollary, we have $\Ind \cE_r(J) = \spkl{\cE_r(J)} \subseteq J \subseteq \Ind \cE_r(J)$.
\eproof

At this point, we remark that associating $\spkl{I_r}$ with an (invariant) ideal $I_r$ of $D_r$ is also a natural way of constructing ideals of $C^*_r(P)$ from those of $D_r$. Indeed, as we will see, this process is to a certain extent even more natural, at least for our purposes. But first, we observe that the assignment $I_r \to \spkl{I_r}$ is also one-to-one (under the condition that $I_r$ is invariant):
\blemma
Given an invariant ideal $I_r$ of $D_r$, we have $\spkl{I_r} \cap D_r = I_r$.
\elemma
\bproof
As we always have $\spkl{I_r} \subseteq \Ind I_r$, our claim follows from $I_r \subseteq \spkl{I_r} \cap D_r \subseteq (\Ind I_r) \cap D_r = I_r$.
\eproof

By our assumptions that $\cJ$ is independent and that $\ping$ is Toeplitz, we know that $C^*_r(P)$ is isomorphic to the full corner $E_P(D_P^G \rtimes_{\tau,r} G)E_P$ of the reduced crossed product $D_P^G \rtimes_{\tau,r} G$. Thus there is a one-to-one correspondence between ideals of $C^*_r(P)$ and ideals of $D_P^G \rtimes_{\tau,r} G$ given by $D_P^G \rtimes_{\tau,r} G \triangleright J \ma J \vert_P \triangleleft C^*_r(P)$. Here $J \vert_P$ is the ideal of $C^*_r(P)$ which corresponds to $E_P J E_P$ under the canonical identification $C^*_r(P) \cong E_P(D_P^G \rtimes_{\tau,r} G)E_P$ provided by Corollary~\ref{full-red}. But even more, we also know that $J \vert_P$ is again isomorphic to a full corner of $J$, namely $E_P J E_P$.

Given an invariant ideal $I_r$ of $D_r$, our present goal is to find a $G$-invariant ideal $I_P^G$ of $D_P^G$ such that $(I_P^G \rtimes_{\tau,r} G) \vert_P = \spkl{I_r}$. A natural candidate for $I_P^G$ would be the smallest $G$-invariant ideal of $D_P^G$ which contains $I_r$.
\bdefin
We set
\bgloz
  I_P^G \defeq \clspan(\menge{\tau_{g_1}(x_1) \dotsm \tau_{g_n}(x_n) \cdot d}{n \in \Zz_{\geq 1}, g_i \in G, x_i \in I_r, d \in D_P^G}.
\egloz
\edefin
We observe that it is an easy consequence of the construction of $I_P^G$ that in $\Omega_P^G$, we have $\Spec(I_P^G) = (\Spec I_r) \cdot G$. Here and in the sequel, we identify $\Omega$ with a subspace of $\Omega_P^G$ via the map $c^*$ from \S~\ref{groupoids}.

\blemma
\label{IP=IPG}
If $\ping$ is Toeplitz, then the following hold:
\begin{enumerate}
\item[(i)] $(I_P^G \rtimes_{\tau,r} G) \vert_P = \spkl{I_r}$,
\item[(ii)] $E_P I_P^G E_P = I_r$,
\item[(iii)] for all $g \in G$, $E_P \tau_g(I_r) E_P \subseteq I_r$,
\item[(iv)] $\Spec I_r = \Spec I_P^G \cap \Omega$,
\item[(v)] $\Omega_P^G \setminus \Spec I_P^G = (\Omega \setminus \Spec I_r) \cdot G$.
\end{enumerate}
\elemma
\bproof
We first prove that these conditions are equivalent if $\ping$ is Toeplitz, and then we show that the Toeplitz condition for $\ping$ implies (iii).

To see \an{(i) $\Rarr$ (ii)}, note that $\spkl{I_r} \cap D_r = I_r$ in $C^*_r(P)$ implies that we have $\spkl{I_r}_{E_P(D_P^G \rtimes_{\tau,r} G)E_P} \cap D_r = I_r$ in $E_P(D_P^G \rtimes_{\tau,r} G)E_P$. Thus if (i) holds, i.e. if $E_P(I_P^G \rtimes_{\tau,r} G)E_P = \spkl{I_r}_{E_P(D_P^G \rtimes_{\tau,r} G)E_P}$, then $I_r \subseteq E_P I_P^G E_P \subseteq E_P(I_P^G \rtimes_{\tau,r} G)E_P \cap D_r = I_r$.

\an{(ii) $\Rarr$ (iii)} is clear as $\tau_g(I_r) \subseteq I_P^G$.

To prove \an{(iii) $\Rarr$ (i)}, we first observe that $(I_P^G \rtimes_{\tau,r} G) \vert_P \supseteq \spkl{I_r}$ always holds as $I_P^G \supseteq I_r$. It remains to prove that (iii) implies the reverse inclusion. Upon identifying $C^*_r(P)$ with $E_P(D_P^G \rtimes_{\tau,r} G)E_P$, we have to prove, assuming (iii), that $E_P(I_P^G \rtimes_{\tau,r} G)E_P \subseteq \spkl{I_r}_{E_P(D_P^G \rtimes_{\tau,r} G)E_P}$. Take a generator of $I_P^G$, say $\tau_{g_1}(x_1) \dotsm \tau_{g_n}(x_n) \cdot d$. Then for all $g \in G$, $E_P \tau_{g_1}(x_1) \dotsm \tau_{g_n}(x_n) \cdot d \cdot U_g E_P = (E_P \tau_{g_1}(x_1) E_P) \cdot (E_P \tau_{g_2}(x_2) \dotsm \tau_{g_n}(x_n) \cdot d \cdot U_g E_P)$, and since $E_P \tau_{g_1}(x_1) E_P$ lies in $I_r$ by (iii), we conclude $E_P \tau_{g_1}(x_1) \dotsm \tau_{g_n}(x_n) \cdot d \cdot U_g E_P \in \spkl{I_r}_{E_P(D_P^G \rtimes_{\tau,r} G)E_P}$.

\an{(iii) $\Rarr$ (iv)}: We always have \an{$\subseteq$}. To prove \an{$\supseteq$}, take $\chi \in \Omega$ such that $\chi \vert_{I_P^G} \neq 0$. Then we can find $g \in G$ and $x \in I_r$ with $\chi(\tau_g(x)) \neq 0$. Thus $\chi(E_P \tau_g(x) E_P) \neq 0$. As $E_P \tau_g(x) E_P$ lies in $I_r$ by (iii), we conclude $\chi \in \Spec I_r$.

\an{(iv) $\Rarr$ (v)}: The inclusion \an{$\subseteq$} is easy to see. The other one (\an{$\supseteq$}) follows from (iv) and $G$-invariance of $\Omega_P^G \setminus \Spec I_P^G$.

\an{(v) $\Rarr$ (iii)}: We have to show that whenever $\chi \in \Omega_P^G$ satisfies $\chi \vert_{I_r} = 0$, then for all $g \in G$, $\chi \vert_{E_p \tau_g(I_r) E_P} = 0$ must hold as well. Take $\chi$ such that $\chi \vert_{I_r} = 0$. If $\chi(E_P) = 0$, there is nothing to show. Hence we may assume $\chi(E_P) = 1$, i.e. $\chi \in \Omega$. This means that $\chi \in \Omega \setminus \Spec I_r$. By (v), we conclude that $\chi \notin \Spec(I_P^G)$.

It remains to prove that the Toeplitz condition for $\ping$ implies (iii). Given $g \in G$ with $E_P \lambda_g E_P \neq 0$, the Toeplitz condition yields $p_1, q_1, \dotsc, p_n, q_n \in P$ such that $E_P \lambda_g E_P = V_{p_1}^* V_{q_1} \dotsm V_{p_n}^* V_{q_n}$. Then, for every invariant ideal $I_r$ of $D_r$, we have
\bgloz
  E_P \tau_g(I_r) E_P = E_P \lambda_g E_P I_r E_P \lambda_g^* E_P = \Ad(V_{p_1}^* V_{q_1} \dotsm V_{p_n}^* V_{q_n})(I_r) \subseteq I_r
\egloz
as $I_r$ is invariant. The case $E_P \lambda_g E_P = 0$ is trivial.
\eproof

\bremark
In particular, if $\ping$ is Toeplitz, then $(\Omega \setminus \Spec I_r) \cdot G$ is closed in $\Omega_P^G$.
\eremark

Now, the same arguments used in the proof of Theorem~\ref{thm1} and Theorem~\ref{thm2} give us the following
\bprop
\label{inv-id-full}
Assume that $\cJ$ is independent and that $\ping$ is Toeplitz. Let $I_r$ be an invariant ideal of $D_r$, and let $I$ be the corresponding ideal of $D$ such that $\lambda(I) = I_r$. Then the ideal $\spkl{I_r}$ of $C^*_r(P)$ generated by $I_r$ is isomorphic to the full corner of $I_P^G \rtimes_{\tau,r} G$ determined by the characteristic function $\1z_{\Spec I_r}$ of $\Spec I_r \subseteq \Spec (I_P^G)$, and the ideal $\spkl{I}$ of $C^*_s(P)$ generated by $I$ is isomorphic to the full corner of $I_P^G \rtimes_{\tau} G$ given by $\1z_{\Spec I_r}$.

Moreover, the following are equivalent: 
\begin{enumerate}
\item[(i$_I$)] $\spkl{I}_{C^*_s(P)}$ is nuclear,
\item[(ii$_I$)] $\spkl{I_r}_{C^*_r(P)}$ is nuclear,
\item[(iii$_I$)] the transformation groupoid $\Spec (I_P^G) \rtimes G = ((\Spec I_r) \cdot G) \rtimes G$ is amenable.
\end{enumerate}

Either of these conditions implies that $\lambda: \spkl{I}_{C^*_s(P)} \to \spkl{I_r}_{C^*_r(P)}$ is faithful.

For the corresponding quotients, we have that $C^*_s(P) / \spkl{I}$ is isomorphic to the full corner of $(D_P^G / I_P^G) \rtimes_{\tau} G \cong C_0(\Omega_P^G \setminus \Spec(I_P^G)) \rtimes_{\tau} G$ determined by the characteristic function $\1z_{\Omega \setminus \Spec I_r}$ of $\Omega \setminus \Spec I_r \subseteq \Omega_P^G \setminus \Spec(I_P^G)$. Moreover, if the sequence $0 \to C_0(\Spec(I_P^G)) \rtimes_{\tau,r} G \to C_0(\Omega_P^G) \rtimes_{\tau,r} G \to C_0(\Omega_P^G \setminus \Spec(I_P^G)) \rtimes_{\tau,r} G \to 0$ is exact, then also $C^*_r(P) / \spkl{I_r}$ is isomorphic to the full corner of $(D_P^G / I_P^G) \rtimes_{\tau,r} G \cong C_0(\Omega_P^G \setminus \Spec(I_P^G)) \rtimes_{\tau,r} G$ associated with $\1z_{\Omega \setminus \Spec I_r}$, and the following are equivalent:
\begin{enumerate}
\item[(i$_Q$)] $C^*_s(P) / \spkl{I}$ is nuclear,
\item[(ii$_Q$)] $C^*_r(P) / \spkl{I_r}$ is nuclear,
\item[(iii$_Q$)] the transformation groupoid $(\Omega_P^G \setminus \Spec(I_P^G)) \rtimes G = ((\Omega \setminus (\Spec I_r)) \cdot G) \rtimes G$ is amenable;
\end{enumerate}
and either of these conditions implies that $\lambda: C^*_s(P) / \spkl{I} \to C^*_r(P) / \spkl{I_r}$ is faithful.
\eprop

We also mention the following useful consequence:
\blemma
\label{invid=invid}
If $\ping$ satisfies the Toeplitz condition, then the maps $I_r \ma I_P^G$ and $E_P J E_P \mafr J$ are mutually inverse, inclusion-preserving bijections between the sets of invariant ideals of $D_r$ and $G$-invariant ideals of $D_P^G$.
\elemma
\bproof
By the Toeplitz condition, we have $E_P I_P^G E_P = I_r$. To check $(E_P J E_P)_P^G = J$, note that $E_P(J \rtimes_{\tau,r} G)E_P = \spkl{E_P J E_P} = E_P((E_P J E_P)_P^G \rtimes_{\tau,r} G)E_P$. As $E_P$ is a full projection in $D_P^G \rtimes_{\tau,r} G$, we conclude that $J \rtimes_{\tau,r} G = (E_P J E_P)_P^G \rtimes_{\tau,r} G$, hence $J = (E_P J E_P)_P^G$.
\eproof

\subsection{Invariant spectral subsets}
\label{des-spec}

As ideals of $D_r$ correspond to subsets of $\Omega = \Spec D_r$, we now describe $\Omega$ explicitly in terms of the family of constructible ideals $\cJ$, and we also describe the action of $P$. This is just an application of our observations in \S~\ref{fam-sub} because of our standing assumption that $\cJ$ is independent.

Let $\Sigma$ be the set of non-empty $\cJ$-valued filters as introduced before Corollary~\ref{Jvf}, equipped with the topology introduced after Corollary~\ref{Jvf}.

\blemma
We can identify $\Omega$ with $\Sigma$ via $\omega: \Omega \to \Sigma$, $\chi \ma \menge{X \in \cJ}{\chi(E_X)=1}$.

For all $p \in P$, the map $\Sigma \to \Sigma$, $\cF \ma \menge{X}{p^{-1} X \in \cF} \eqdef p \cF$ gives rise to a homeomorphism $\Sigma \cong p \Sigma = \menge{\cF \in \Sigma}{pP \in \cF}$. Let $p^{-1}: p \Sigma \to \Sigma, p \cF \ma \cF$ be its inverse. Define $\sigma_p: C(\Sigma) \to C(\Sigma)$, $\sigma_p(d)(\cF) = d(p^{-1} \cF)$ if $\cF$ lies in $p \Sigma$ and $\sigma_p(d)(\cF) = 0$ if $\cF$ does not lie in $p \Sigma$; and $\sigma_{p^*}: C(\Sigma) \to C(\Sigma)$, $\sigma_{p^*}(d)(\cF) = d(p \cF)$.

Then the homeomorphism $\omega: \Omega \to \Sigma$ induces an identification $\omega^*: C(\Sigma) \to D_r$ such that for every $p \in P$, the diagrams
\bgloz
  \begin{CD}
  C(\Sigma) @>>> D_r \\
  @V \sigma_p VV @VV {\Ad(V_p)} V \\
  C(\Sigma) @>>> D_r
  \end{CD}
\egloz
and
\bgloz
  \begin{CD}
  C(\Sigma) @>>> D_r \\
  @V \sigma_{p^*} VV @VV {\Ad(V_p^*)} V \\
  C(\Sigma) @>>> D_r
  \end{CD}
\egloz
commute.
\elemma
\bproof
This is straightforward to check.
\eproof

\bcor
\label{invid-invs}
An ideal $I_r$ of $D_r$ is invariant if and only if for all $p \in P$, we have $p \omega(\Spec I_r) \subseteq \omega(\Spec I_r)$ and $p^{-1}(\omega(\Spec I_r) \cap p \Sigma) \subseteq \omega(\Spec I_r)$.
\ecor

\bdefin
A subset $C$ of $\Sigma$ is called invariant if for all $p \in P$, the conditions $pC \subseteq C$ and $p^{-1}(C \cap p\Sigma) \subseteq C$ are satisfied.
\edefin

\subsection{The boundary action}
\label{bd-act}

Finally, let us have a look at the boundary of $\Omega$. Recall the definition of the boundary $\partial \Omega$ from \S~\ref{fam-sub}.

\bdefin
Let $\Sigma_{\max}$ be the set of all $\cJ$-valued ultrafilters, and let $\partial \Sigma$ be the closure $\overline{\Sigma_{\max}}$ of $\Sigma_{\max}$ in $\Sigma$. We set $\Omega_{\max} = \omega^{-1}(\Sigma_{\max})$ and $\partial \Omega = \omega^{-1}(\partial \Sigma)$.
\edefin

We note that this definition is essentially the one from \cite{La1}, extended from the case of quasi-lattice ordered groups to our situation.

\blemma
\label{bound-min}
$\partial \Sigma$ is the minimal non-empty closed invariant subset of $\Sigma$.
\elemma
\bproof
Choose $\cF \in \Sigma_{\max}$ and take $p \in P$.

We claim $p \cF \in \Sigma_{\max}$. Assume that there exists $\cF' \in \Sigma$ such that $p \cF \subseteq \cF'$. Then $pP \in \cF'$ so that $\cF' \in p\Sigma$. Thus $p^{-1} \cF'$ is an element of $\Sigma$ such that $\cF = p^{-1} p \cF \subseteq p^{-1} \cF'$. As $\cF$ is maximal, this implies $\cF = p^{-1} \cF'$. Hence $p \cF = p p^{-1} \cF' = \cF'$.

Next we claim that $p^{-1} \cF \in \Sigma_{\max}$ if $\cF$ lies in $p \Sigma$. If $p^{-1} \cF \subseteq \cF'$ for some $\cF' \in \Sigma$, then $\cF = p p^{-1} \cF \subseteq p \cF'$ implies $\cF = p \cF'$ and thus $p^{-1} \cF = p^{-1} p \cF' = \cF'$.

Thus we have seen that $\Sigma_{\max}$ is invariant. As $\partial \Sigma$ is the closure of $\Sigma_{\max}$, we conclude that $p (\partial \Sigma) \subseteq \partial \Sigma$ for all $p \in P$. As we know that $p \Sigma$ is clopen in $\Sigma$, we also deduce $p^{-1}(\partial \Sigma \cap p \Sigma) \subseteq \partial \Sigma$. Therefore $\partial \Sigma$ is invariant.

To prove minimality, let $\emptyset \neq C$ be a closed invariant subset of $\Sigma$. Take $\cF \in C$ arbitrary, and choose some $\cF_{\max} \in \Sigma_{\max}$. For every $X \in \cF_{\max}$, choose $x \in X$ ($X \neq \emptyset$). Then $xP \in x\cF$ implies that $X$ lies in $x\cF$ as $xP \subseteq X$ ($X$ is a right ideal). Set $\cF_X \defeq x\cF$. Ordering elements in $\cF_{\max}$ by inclusion (i.e. we set $X_1 \geq X_2$ if $X_1 \subseteq X_2$), we obtain a net $(\cF_X)_{X \in \cF_{\max}}$ in $C$. As $\Sigma$ is compact ($\Omega$ is compact), we may assume, after passing to a convergent subnet if necessary, that $(\cF_X)_X$ converges to an element $\cF_m$ of $\Sigma$. As $C$ is closed, $\cF_m$ must lie in $C$. Moreover, for every $X \in \cF_{\max}$, we have that $X' \geq X$ implies $X \in \cF_{X'}$. Thus $X$ lies in $\cF_m$. We conclude $\cF_{\max} \subseteq \cF_m$, hence by maximality, $\cF_{\max} = \cF_m$ lies in $C$. Thus $\partial \Sigma = \overline{\Sigma_{\max}}$ lies in $C$, and we are done.
\eproof

As immediate consequences, we obtain
\bcor
$V_r(\partial \Omega) = \menge{d \in D_r}{\chi(d) = 0 \fa \chi \in \partial \Omega}$ is the maximal invariant proper ideal of $D_r$.
\ecor
\bproof
This is a direct consequence of minimality of $\partial \Omega$ and Corollary~\ref{invid-invs}.
\eproof
\bcor
\label{min-act}
If $\ping$ satisfies the Toeplitz condition, then the $G$-action on $(\partial \Omega) \cdot G$ is minimal.
\ecor
\bproof
This is a direct consequence of the previous corollary and Lemma~\ref{invid=invid}.
\eproof

Moreover, we deduce the analogue of \cite{La1}, Proposition~4.3:
\bcor
Given a proper ideal $J$ of $C^*_r(P)$ such that $\cE_r(J) \subseteq J$, we always have $J \subseteq \Ind V_r(\partial \Omega)$.
\ecor
\bproof
As $J$ is a proper ideal of $C^*_r(P)$ and since $\cE_r(J) \subseteq J$, $\cE_r(J)$ is an invariant proper ideal of $D_r$. Thus $\cE_r(J) \subseteq V_r(\partial \Omega)$. Moreover, for every $x \in J$, we have $\cE_r(x^*x) \in \cE_r(J)$. Hence every element in $J$ lies in $\menge{x \in C^*_r(P)}{\cE_r(x^*x) \in V_r(\partial \Omega)} = \Ind V_r(\partial \Omega)$.
\eproof

Using similar ideas as in \cite{La1}, we now investigate when the action of $G$ on $(\partial \Omega) \cdot G \subseteq \Omega_P^G$ is topologically free and a local boundary action (the first notion is introduced in \cite{Ar-Sp}, and the second one is introduced in \cite{La-Sp}).

First of all, we set
\bgloz
  G_0 \defeq \menge{g \in G}{(g \cdot P) \cap X \neq \emptyset \text{ and } (g^{-1} \cdot P) \cap X \neq \emptyset \fa \emptyset \neq X \in \cJ}.
\egloz
Clearly, $G_0 = \menge{g \in G}{(g \cdot P) \cap (pP) \neq \emptyset \text{ and } (g^{-1} \cdot P) \cap (pP) \neq \emptyset \fa p \in P}$. Moreover, we have
\blemma
$G_0$ is a subgroup of $G$.
\elemma
\bproof
Take $g_1$, $g_2$ in $G_0$. Then for all $\emptyset \neq X \in \cJ$, we have $((g_1 g_2) \cdot P) \cap X = g_1 \cdot ((g_2 \cdot P) \cap (g_1^{-1} \cdot X)) \supseteq g_1 \cdot ((g_2 \cdot P) \cap (g_1^{-1} \cdot X)) \cap (g_1 \cdot P) = g_1 \cdot ((g_2 \cdot P) \cap ((g_1^{-1} \cdot X) \cap P))$. Now $(g_1^{-1} \cdot X) \cap P = g_1^{-1} \cdot (X \cap (g_1 \cdot P)) \neq \emptyset$. Thus there exists $x \in P$ such that $x \in (g_1^{-1} \cdot X) \cap P$. Hence $xP \subseteq (g_1^{-1} \cdot X) \cap P$. Thus $\emptyset \neq g_1 \cdot ((g_2 \cdot P) \cap (xP)) \subseteq ((g_1 g_2) \cdot P) \cap X$.
\eproof

\bprop
\label{top-free}
$G$ acts topologically freely on $(\partial \Omega) \cdot G$ if and only if $G_0$ acts topologically freely on $(\partial \Omega) \cdot G$.
\eprop
\bproof
\an{$\Rarr$} is clear. For \an{$\Larr$}, assume that for $e \neq g \in G$, we have that the fix point set of $(\partial \Omega) \cdot G$ under $g$ has non-empty interior, i.e. $\eckl{\text{Fix}_{(\partial \Omega) \cdot G}(g)}^{\circ} \neq \emptyset$.  Thus there exists an open subset $U$ of $(\partial \Omega) \cdot G$ such that $U \subseteq \text{Fix}_{(\partial \Omega) \cdot G}(g)$. As $\partial \Omega = \overline{\Omega_{\max}}$, we deduce $(\partial \Omega) \cdot G \subseteq \overline{\Omega_{\max} \cdot G}$. Therefore there exists $\chi \in (\Omega_{\max} \cdot G) \cap U$. Choose $h \in G$ such that $\chi h$ lies in $\Omega_{\max}$.

Now assume that $G_0$ acts topologically freely on $(\partial \Omega) \cdot G$. Then for all $x \in G$, $x^{-1} g x$ cannot lie in $G_0$ as $\text{Fix}_{(\partial \Omega) \cdot G}(x^{-1}gx) = \text{Fix}_{(\partial \Omega) \cdot G}(g) \cdot x^{-1}$. Now take any $X \in \cJ$ such that $(\chi h)(E_X) = 1$. Moreover, let $x \in X$. Then $x^{-1} h^{-1} g h x$ does not lie in $G_0$. Thus there exists $p \in P$ such that $((x^{-1} h^{-1} g h x) \cdot P) \cap (pP) = \emptyset$ or $((x^{-1} h^{-1} g^{-1} h x) \cdot P) \cap (pP) = \emptyset$. In either case, we choose $\chi_X \in \Omega_{\max}$ such that $\chi_X(E_{xpP})=1$. If $((x^{-1} h^{-1} g h x) \cdot P) \cap (pP) = \emptyset$, then $(xpP) \cap ((h^{-1} g h x) \cdot P)) = \emptyset$. This implies $\chi_X(E_{(h^{-1} g h x) \cdot P}) = 0$. Thus $\chi_X \cdot (h^{-1} g h) \neq \chi_X$. Similarly, we obtain from $((x^{-1} h^{-1} g^{-1} h x) \cdot P) \cap (pP) = \emptyset$ that $\chi_X \cdot (h^{-1} g^{-1} h) \neq \chi_X$. In any case, we obtain $\chi_X \cdot (h^{-1} g) \neq \chi_X h^{-1}$, hence $\chi_X h^{-1}$ does not lie in $\text{Fix}_{(\partial \Omega) \cdot G}(g)$.

As against that, we have found for all $X \in \omega(\chi h)$ a character $\chi_X \in \Omega_{\max}$ with $\chi_X(E_X) = 1$. Thus ordering $X \in \omega(\chi h)$ by inclusion as in the proof of Lemma~\ref{bound-min}, we obtain a net $(\chi_X)_X$ in $\Omega_{\max} \subseteq \Omega$. By passing over to a convergent subnet, we may assume that $\lim_X \chi_X = \ti{\chi} \in \Omega$. Hence $\ti{\chi}(E_X)=1$ for all $X \in \omega(\chi h)$. This implies $\omega(\chi h) \subseteq \omega(\ti{\chi})$, hence $\ti{\chi} = \chi h$ as $\chi h$ lies in $\Omega_{\max}$. The conclusion is that $\lim_X \chi_X h^{-1} = \chi$. But we have seen $\chi_X h^{-1} \notin \text{Fix}_{(\partial \Omega) \cdot G}(g)$, and we also know $\chi \in \eckl{\text{Fix}_{(\partial \Omega) \cdot G}(g)}^{\circ}$. This is a contradiction.
\eproof

\bprop
\label{loc-bd-act}
If $P$ is not left reversible, then the $G$-action on $(\partial \Omega) \cdot G$ is a local boundary action in the sense of \cite{La-Sp}.
\eprop
\bproof
We have to show that for every non-empty open subset $U$ of $(\partial \Omega) \cdot G$, there exists an open subset $\Delta \subseteq U$ and an element $g \in G$ such that $\overline{\Delta} g \subsetneq \Delta$.

Let $U$ be as above. As $\overline{\Omega_{\max}} = \partial \Omega$, we can find $\chi \in \Omega_{\max}$ and $h \in G$ such that $\chi h \in U$, i.e. $\chi \in (U h^{-1}) \cap \Omega$. As $\Omega$ is open in $\Omega_P^G$, we can find $X$ in $\cJ$ and $X_1$, ..., $X_n$ in $\cJ$ such that $V = \menge{\psi \in (\partial \Omega) \cdot G}{\psi(E_X)=1, \psi(E_{X_i})=0 \fa 1 \leq i \leq n}$ is contained in $U h^{-1}$ and that $\chi \in V$ (see \eqref{basis-top}). The latter condition means that $\chi(E_X)=1$ and $\chi(E_{X_i})=0$ for all $1 \leq i \leq n$. As $\chi$ lies in $\Omega_{\max}$, we conclude that for all $1 \leq i \leq n$, there must be $X'_i$ in $\cJ$ such that $\chi(E_{X'_i})=1$ and $X_i \cap X'_i = \emptyset$ (see Remark~\ref{F-in-Sigma}). Thus setting $\ti{X} = X \cap (\bigcap_{i=1}^n X'_i) \neq \emptyset$, we see that for any $\psi \in (\partial \Omega) \cdot G$, $\psi(E_{\ti{X}})=1$ implies $\psi \in U h^{-1}$. Choose $x \in \ti{X}$. As $P$ is not left reversible, we can find $p$ and $q$ in $P$ such that $(pP) \cap (qP) = \emptyset$. Now set $\Delta' = \menge{\psi \in (\partial \Omega) \cdot G}{\psi(E_{xP})=1}$. $\Delta'$ is clopen in $(\partial \Omega) \cdot G$. As $xP \subseteq \ti{X}$, we conclude that $\Delta' \subseteq U h^{-1}$. Set $g' = xp^{-1}x^{-1} \in G$. Then $\psi' \in \Delta' g'$ implies that there exists $\psi \in (\partial \Omega) \cdot G$ such that $\psi(E_{xP})=1$ and $\psi' = \psi g'$. Thus $\psi' \in (\partial \Omega) \cdot G$ and $\psi'(E_{xpP}) = \psi(E_{g' \cdot (xpP)}) = \psi(E_{xP}) = 1$. Thus $\psi'(E_{xP})=1$ as $xpP \subseteq xP$. Hence $\psi'$ lies in $\Delta'$. This shows $\Delta' g' \subseteq \Delta'$. But we can now choose $\psi' \in \Delta'$ with $\psi'(E_{xqP})=1$. If $\psi'$ lies in $\Delta' g'$, then there exists $\psi \in \Delta$ such that $\psi'=\psi g'$. Then $\psi'(E_{xpP}) = \psi'(E_{g' \cdot (xpP)}) = \psi'(E_{xP}) = 1$. But this contradicts $(xqP) \cap (xpP) = x((qP) \cap (pP)) = \emptyset$. Hence $\psi'$ does not lie in $\Delta' g'$, and we have proven $\Delta' g' \subsetneq \Delta'$. Setting $\Delta \defeq \Delta' h$ and $g = h^{-1} g' h$, we are done.
\eproof

\bremark
Proposition~\ref{loc-bd-act} clarifies the final remark in \cite{La1}, where it is pointed out that the boundary action should be \an{a boundary action in the sense generalizing that of \cite{La-Sp}}.
\eremark

\bcor
\label{pisun}
Assume that $G$ is countable, that $\cJ$ is independent and that $\ping$ satisfies the Toeplitz condition. Let $V(\partial \Omega)$ be the ideal of $D$ such that $\lambda(V(\partial \Omega)) = V_r(\partial \Omega)$. 

The boundary quotient $C^*_s(P) / \spkl{V(\partial \Omega)}$ is a unital UCT Kirchberg algebra if and only if the following hold:
\begin{itemize}
\item $P \neq \gekl{e}$, 
\item $G$ acts amenably on $(\partial \Omega) \cdot G$,
\item $G_0$ acts topologically freely on $(\partial \Omega) \cdot G$.
\end{itemize}
\ecor
\bproof
As $G$ is countable, so is $P$. Thus $C_0((\partial \Omega) \cdot G) \rtimes_{\tau} G$ is separable. By Proposition~\ref{inv-id-full}, $C^*_s(P) / \spkl{V(\partial \Omega)}$ and $C_0((\partial \Omega) \cdot G) \rtimes_{\tau} G$ are stably isomorphic. Hence $C^*_s(P) / \spkl{V(\partial \Omega)}$ is a Kirchberg algebra if and only if $C_0((\partial \Omega) \cdot G) \rtimes_{\tau} G$ is a Kirchberg algebra. By \cite{Br-Oz}, Chapter~5, Theorem~6.18 and the last corollary of \cite{Ar-Sp}, $C_0((\partial \Omega) \cdot G) \rtimes_{\tau} G$ is nuclear and simple if and only if $G$ acts on $(\partial \Omega) \cdot G$ amenably and topologically freely. Here we have used Corollary~\ref{min-act} which tells us that the $G$-action on $(\partial \Omega) \cdot G$ is minimal. Topological freeness of the $G$-action is equivalent to topological freeness of the $G_0$-action by Proposition~\ref{top-free}. To complete our proof, first observe that clearly, $P$ has to be non-trivial if $C_0((\partial \Omega) \cdot G) \rtimes_{\tau} G$ is purely infinite. This settles the implication \an{$\Rarr$}. For the converse, we show that our assumptions that $P \neq \gekl{e}$ and that $G_0$ acts topologically freely on $(\partial \Omega) \cdot G$ imply that $P$ is not left reversible: If $P$ were left reversible, i.e. if every non-empty $X_1$, $X_2$ in $\cJ$ have non-empty intersection, then $\partial \Omega$ would consist of only one point, namely the $\cJ$-valued ultrafilter $\cJ\reg$ consisting of all non-empty elements in $\cJ$. Also, if $P$ were left reversible, then we would have $P \subseteq G_0$. Since every element in $P$ obviously leaves $\cJ\reg$ fixed, and by our assumption that $P \neq \gekl{e}$, we conclude that $G_0$ cannot act topologically freely on $(\partial \Omega) \cdot G$ if $P$ were left reversible. Hence Proposition~\ref{loc-bd-act} implies that the $G$-action on $(\partial \Omega) \cdot G$ is a local boundary action. With the help of Theorem~9 from \cite{La-Sp}, this settles the reverse direction \an{$\Larr$}.
\eproof

\section{Examples}
\label{ex}

\subsection{Quasi-lattice ordered groups}
\label{qlo}

Recall from \cite{Ni1} that a pair $(G,P)$ consisting of a subsemigroup $P$ of a group $G$ is called quasi-lattice ordered if
\begin{itemize}
\item[(QL0)] $P \cap P^{-1} = \gekl{e}$,
\item[(QL1)] for all $g \in G$, the intersection $P \cap (g \cdot P)$ is either empty or of the form $pP$ for some $p \in P$.
\end{itemize}
As observed in \cite{Cr-La1}, \S~3 (Remark~8), (QL1) implies
\begin{itemize}
\item[(QL2)] For all $p$, $q$ in $P$, the intersection $(pP) \cap (qP)$ is either empty or of the form $rP$ for some $r \in P$.
\end{itemize}
In the sequel, we will most of the time only use (QL1) and (QL2).

First of all, we observe that for every such $\ping$ satisfying (QL2), we have $\cJ = \menge{pP}{p \in P} \cup \gekl{\emptyset}$. In this sense, the ideal structure (or rather the structure of the constructible right ideals) is very simple. It is immediate that $\cJ$ is independent. Moreover, $\ping$ satisfies the Toeplitz condition (compare also \cite{Ni1}, \S~2.4). Namely, take $g \in G$. If $E_P \lambda_g E_P \neq 0$, then there exists $p \in P$ such that $P \cap (g \cdot P) = pP$ by (QL1). Thus there is $q \in P$ with $gq=p$, hence $g=pq^{-1}$. It then follows that $E_P \lambda_g E_P = E_{P \cap (g \cdot P)} \lambda_p \lambda_{q^{-1}} E_P = E_{pP} \lambda_p \lambda_{q^{-1}} E_P = (E_P \lambda_p E_P)(E_P \lambda_{q^{-1}} E_P) = V_p V_q^*$. Therefore, all our results apply.

As mentioned in the introduction, quasi-lattice ordered groups and their semigroup C*-algebras have been studied intensively, for instance in \cite{Ni1}, \cite{Ni2}, \cite{La-Rae}, \cite{La1}, \cite{E-L-Q}, \cite{Cr-La1} and \cite{Cr-La2}. The full semigroup C*-algebras have been described as semigroup crossed products by endomorphisms in \cite{La-Rae}. Moreover, both full and reduced semigroup C*-algebras can be described as partial crossed products of the corresponding groups. This gives yet another description which is not discussed here, but which is certainly closely related to \S~\ref{var-des}. The induced ideals of reduced semigroup C*-algebras have been studied in \cite{Ni1}. The $\Theta_i$ we introduced in Corollary~\ref{Theta} can be viewed as a substitute for the positive definite functions $\theta_i$ introduced in \cite{Ni1}, \S~4.5. And the reader will see that for Proposition~\ref{nuc-ind-id}, we have essentially adapted A. Nica's proof of the proposition in \S~6.1 of \cite{Ni1}. The boundary of the spectrum was introduced in \cite{La1} and studied in \cite{La1}, \cite{Cr-La2}. Our discussion of the boundary action in \S~\ref{bd-act} is modelled after \cite{La1} and \cite{Cr-La2}.

Before we come to an explicit example, let us first show how the analysis in \cite{La-Rae} can be extended. Namely, we obtain a strengthening of Proposition~6.6 in \cite{La-Rae} with essentially the same proof as in \cite{La-Rae}. We point out that the conclusion in this proposition should read \an{If $\cG$ is amenable, then $(G,P)$ is amenable} (compare also Remark~17 in \cite{Cr-La1}). Let us start with the following
\bprop[\cite{La-Rae}, Lemma~4.1 for arbitrary coefficients]
\label{LR1}
Let $(G,P)$ and $(H,Q)$ be quasi-lattice ordered. Assume that $\varphi: G \to H$ is a group homomorphism such that $\varphi(P) \subseteq Q$ and whenever $x$, $y$ in $P$ satisfy $(xP) \cap (yP) \neq \emptyset$, then
\bgln
\label{controll1}
  && \varphi(x) = \varphi(y) \LRarr x=y, \\
\label{controll2}
  && \text{for } z \in P \text{ such that } (xP) \cap (yP) = zP \text{, } (\varphi(x) Q) \cap (\varphi(y) Q) = \varphi(z) Q.
\egln
Moreoever, let $\alpha$ be a $G$-action on a C*-algebra $A$.

Then $B \defeq \clspan(\menge{\iota(a) \overline{v_x v_y^*}}{a \in A; x, y \in P \text{ with } \varphi(x)=\varphi(y)}$ is a sub-C*-algebra of $A \rta_{\alpha,s} P$ such that $\lambda_{(A,P,\alpha)} \vert_B: B \to A \rta_{\alpha,r} P$ is faithful.
\eprop
\bproof
Let $F \subseteq Q$ be a finite subset such that whenever $f_1$, $f_2$ in $F$ and $f_3$ in $Q$ satisfy $(f_1 Q) \cap (f_2 Q) = f_3 Q$, then $f_3$ lies in $F$ as well. The set 
\bgl
\label{set-mult?}
  \menge{\iota(a) \overline{v_x v_y^*}}{a \in A; x, y \in P \text{ with } \varphi(x) = \varphi(y) \in F}
\egl
is obviously *-invariant. Moreover, given $\iota(a_1) \overline{v_{x_1} v_{y_1}^*}$ and $\iota(a_2) \overline{v_{x_2} v_{y_2}^*}$ from this set, let $(y_1 P) \cap (x_2 P) = zP$ with $z = y_1 z_1 = y_2 z_2$ for some $z_1$, $z_2$ in $P$. Then
\bglnoz
  \iota(a_1) \overline{v_{x_1} v_{y_1}^*} \iota(a_2) \overline{v_{x_2} v_{y_2}^*}
  &=& \iota(a_1 \alpha_{x_1 y_1^{-1}}(a_2)) \overline{v_{x_1} v_{y_1}^*} \underbrace{\overline{v_{y_1} v_{y_1}^* v_{x_2} v_{x_2}^*}}_{= \overline{v_z v_z^*}}
  \overline{v_{x_2} v_{y_2}^*} \\
  &=& \iota(a_1 \alpha_{x_1 y_1^{-1}}(a_2)) \overline{v_{x_1 z_1} v_{y_2 z_2}^*}.
\eglnoz
Since $\varphi(x_1 z_1) = \varphi(y_1 z_1) = \varphi(z) = \varphi(x_2 z_2) = \varphi(y_2 z_2)$ lies in $F$ by \eqref{controll2}, we have seen that \eqref{set-mult?} is multiplicatively closed. Hence
\bgloz
  B_F \defeq \clspan(\menge{\iota(a) \overline{v_x v_y^*}}{a \in A; x, y \in P \text{ with } \varphi(x) = \varphi(y) \in F}
\egloz
is a sub-C*-algebra of $A \rta_{\alpha,s} P$. As we can write $B = \overline{\bigcup_F B_F}$, we see that $B$ is a sub-C*-algebra of $A \rta_{\alpha,s} P$. Moreover, it suffices to prove faithfulness of $\lambda_{(A,P,\alpha)}$ on $B_F$ for every $F$. Let us first take $F = \gekl{s}$ and consider the representation $\lambda_{(A,P,\alpha)}$ of $B_{\gekl{s}}$ restricted to $\cH \otimes \ell^2(P \cap \varphi^{-1}(\gekl{s})) \subseteq \cH \otimes \ell^2(P)$. Take $\iota(a) \overline{v_x v_y^*} \in B_{\gekl{s}}$. For $z \in \varphi^{-1}(\gekl{s})$, either $z \notin yP$ which implies $\aalphaP (I_\cH \otimes V_x V_y^*) (\xi \otimes \ve_z) = 0$ for all $\xi \in \cH$, or $z$ lies in $yP$. In the latter case, $(zP) \cap (yP) \neq \emptyset$ implies, since $\varphi(z) = \varphi(y) = s$ that $z=y$. Thus $\aalphaP (I_\cH \otimes V_x V_y^*) (\xi \otimes \ve_z) = \delta_{y,z} (\alpha_x^{-1}(a) \xi) \otimes \ve_x$. This means that we have a commutative diagram
\bgloz
  \begin{CD}
  A \otimes_{\max} \cK(\ell^2(P \cap \varphi^{-1}(\gekl{s}))) @>>> A \otimes_{\min} \cK(\ell^2(P \cap \varphi^{-1}(\gekl{s}))) \\
  @VVV @V \subseteq VV \\
  B_{\gekl{s}} @> \lambda_{(A,P,\alpha)} >> \cL(\cH \otimes \ell^2(P \cap \varphi^{-1}(\gekl{s})))
  \end{CD}
\egloz
where the left vertical arrow sends $A \otimes_{\max} \cK(\ell^2(P \cap \varphi^{-1}(\gekl{s}))) \ni a \otimes e_{x,y}$ to $\alpha_x(a) \overline{v_x v_y^*} \in B_{\gekl{s}}$. The upper horizontal arrow is the canonical homomorphism which is an isomorphism as the algebra of compact operators is nuclear. Thus $\lambda_{(A,P,\alpha)}$ is faithful on $B_{\gekl{s}}$.

To go from $B_{\gekl{s}}$ to $B_F$, just proceed as in the proof of Lemma~4.1 in \cite{La-Rae}.
\eproof

\bcor[Proposition~6.6 of \cite{La-Rae} revisited]
Assume that under the hypothesis of the previous proposition, the group $H$ is amenable. Then for every $(A,G,\alpha)$ with $A$ unital, the canonical homomorphism $\lambda_{(A,P,\alpha)}: A \rta_{\alpha,s} P \to A \rta_{\alpha,r} P$ is an isomorphism.
\ecor
\bproof
By \cite{La-Rae}, Proposition~6.1, we have a coaction $A \rta_{\alpha,s} P \to (A \rta_{\alpha,s} P) \otimes_{\max} C^*(H)$ sending $\iota(a) \overline{v_x}$ to $\iota(a) \overline{v_x} \otimes u_{\varphi(x)}$. Here we used that $A \rta_\alpha P \cong A \rta_{\alpha,s} P$ (see \cite{Li2}, \S~3.1) and the crossed product description of $A \rta_\alpha P$ from \cite{Li2}, Lemma~2.15. Thus, as explained in \cite{La-Rae} after Definition~6.3, there exists a conditional expectation $\Psi_{\delta}: A \rta_{\alpha,s} P \to B$ sending $\iota(a) \overline{v_x v_y^*}$ to $\delta_{\varphi(x),\varphi(y)} \iota(a) \overline{v_x v_y^*}$. And by \cite{La-Rae}, Lemma~6.5, this conditional expectation $\Psi_{\delta}$ is faithful if $H$ is amenable. Now let $\cE_r^{(A,P,\alpha)}: A \rta_{\alpha,r} P \to A \otimes D_r$ be the canonical faithful conditional expectation. Then it is straightforward to see that $\cE_r^{(A,P,\alpha)} \circ \lambda_{(A,P,\alpha)} = \cE_r^{(A,P,\alpha)} \circ (\lambda_{(A,P,\alpha)} \vert_B) \circ \Psi_{\delta}$. As the right hand side is faithful by the previous proposition, $\lambda_{(A,P,\alpha)}$ must be faithful.
\eproof

Combining this with Theorem~\ref{thm2} (see also Remark~\ref{unital-suff}), and using Proposition~19 of \cite{Cr-La1}, we obtain
\bcor
Let $(G,P)$ be a quasi-lattice ordered group which admits a map $\varphi$ as in Proposition~\ref{LR1} such that $H$ is amenable. Then $C^*_s(P)$ ($\cong C^*_r(P)$) is nuclear.

In particular, if $(G,P)$ is the graph product of a family of quasi-lattice orders whose underlying groups are amenable, then $C^*_s(P)$ ($\cong C^*_r(P)$) is nuclear.
\ecor

\subsection{Yet another description of Cuntz algebras}
\label{O}

To give an explicit example, consider for $n \geq 2$ the semigroup $\Nz_0^{*n}$, the $n$-fold free product of the natural numbers. Let $p_1$, ..., $p_n$ be the canonical generators of $\Nz_0^{*n}$. The semigroup $\Nz_0^{*n}$ sits inside the free group $\Fz_n$ in a canonical way. This is an example of a quasi-lattice ordered group. It is due to \cite{Ni1}.

Let us now describe $\Omega$ and $\partial \Omega$. As $\cJ = \menge{p \Nz_0^{*n}}{p \in \Nz_0^{*n}} \cup \gekl{\emptyset}$, $\Omega$ can be identified with the set of all, finite or infinite, (reduced) words in the generators $p_1$, ..., $p_n$. Note that we do not allow inverses of the $p_i$ in these words. The topology is the usual restricted product topology. Moreover, the boundary $\partial \Omega$ is precisely the closed subset of all infinite words. $\Omega \setminus (\partial \Omega)$ is then given by the open subset of all finite words. The semigroup $\Nz_0^{*n}$ acts by shifting from the left. The corresponding group action of $\Fz_n$ is given as follows: $\Omega \cdot \Fz_n$ is given by $\Fz_n \cup (\partial \Fz_n)_+$ where $(\partial \Fz_n)_+$ is the set of all infinite words which in reduced form only contain finitely many inverses of the generators $p_1$, ..., $p_n$. The topology is obtained by restricting the canonical one from $\Fz_n \cup (\partial \Fz_n)$. The free group $\Fz_n$ acts by left translations. Moreover, the boundary $(\partial \Omega) \cdot \Fz_n$ is given by $(\partial \Fz_n)_+$.

Let us now turn to the corresponding C*-algebras. Since $\Fz_n$ acts amenably on $\Omega \cdot \Fz_n = \Fz_n \cup (\partial \Fz_n)_+$ (this can be proven for instance as in \cite{Br-Oz}, Chapter~5, \S~1), we do not have to distinguish between full and reduced versions. From the definition, it is clear that $C^*(\Nz_0^{*n})$ is the universal C*-algebra generated by $n$ isometries $v_1$, ..., $v_n$ whose range projections are orthogonal. Therefore $C^*(\Nz_0^{*n})$ is nothing else but the canonical extension of the Cuntz algebra $\cO_n$. Moreover, it is not difficult to see that $\Ind V(\partial \Omega) = \spkl{V(\partial \Omega)}$ is the ideal of $C^*(\Nz_0^{*n})$ generated by the defect projection $1 - \sum_{i=1}^n v_i v_i^*$. Therefore the boundary quotient $C^*(\Nz_0^{*n}) / \spkl{V(\partial \Omega)}$ is canonically isomorphic to $\cO_n$. Passing over to the group crossed products, we obtain
\bglnoz
  && C^*(\Nz_0^{*n}) \sim_M C_0(\Fz_n \cup (\partial \Fz_n)_+) \rtimes \Fz_n, \\
  && \spkl{V(\partial \Omega)} \sim_M C_0(\Fz_n) \rtimes \Fz_n \cong \cK(\ell^2(\Fz_n)), \\
  && C^*(\Nz_0^{*n}) / \spkl{V(\partial \Omega)} \sim_M C_0((\partial \Fz_n)_+) \rtimes \Fz_n.
\eglnoz
The last line gives a description of $\cO_n$ as an ordinary group crossed product by $\Fz_n$ up to Morita equivalence.

Moreover, the group $G_0$ from Proposition~\ref{top-free} is the trivial group in this particular case. Hence Corollary~\ref{pisun} says that $C^*(\Nz_0^{*n}) / \spkl{V(\partial \Omega)}$ is a (unital) UCT Kirchberg algebra. Of course, since we have already observed $C^*(\Nz_0^{*n}) / \spkl{V(\partial \Omega)} \cong \cO_n$, this is not surprising. The point we would like to make is that we did not use anything we already knew about $\cO_n$ to prove all this. So in a way, we have obtained an independent proof of the fact that $\cO_n$ is a UCT Kirchberg algebra (though one has to admit that the proof of pure infiniteness in \cite{La-Sp} is really just the original argument of J. Cuntz).

A similar analysis for the free product $\Nz_0^{* \infty}$ of countably infinitely many copies of the natural numbers yields that $C^*(\Nz_0^{* \infty}) \cong \cO_{\infty}$ is a UCT Kirchberg algebra. In this case, the boundary is everything (i.e. $\Omega = \partial \Omega$) as observed in Remark~3.9 of \cite{La1}. 

\subsection{Left Ore semigroups}
\label{Ore}

Another class of examples is given by left Ore semigroups. Recall that a semigroup $P$ is left Ore if and only if it can be embedded into a group $G$ such that $G = P^{-1} P$. For a left Ore semigroup $P$ with enveloping group $G = P^{-1} P$, $\ping$ always satisfies the Toeplitz condition. Namely, take $g \in G$ and write $g = p^{-1}q$ for $p$, $q$ in $P$. Then $E_P \lambda_g E_P = E_P \lambda_{p^{-1}} \lambda_q E_P = (E_P \lambda_{p^{-1}} E_P) (E_P \lambda_q E_P) = V_p^* V_q$. However, it is not clear whether $\cJ$ is always independent. So we have to assume this.

We remark that setting $\cJ' \defeq \menge{\cap_{i=1}^n p_iP}{n \in \Zz_{\geq 1}, p_i \in P} \cup \gekl{\emptyset}$, we have $\cJ = \menge{q^{-1} X}{q \in P, X \in \cJ'}$. Thus independence of $\cJ$ is equivalent to independence of $\cJ'$. Moreover, in the construction of full semigroup C*-algebras, it actually suffices to consider $\cJ'$ instead of $\cJ$. This is why in \cite{Li1}, only this smaller family $\cJ'$ of right ideals is considered.

Concrete examples of left Ore semigroups are for instance listed in \cite{La2}. Let us briefly discuss the case of $ax+b$-semigroups. Given an integral domain $R \neq \gekl{0}$, we form the semidirect product $R \rtimes R\reg$, where $R\reg = R \setminus \gekl{0}$ acts on the additive group $(R,+)$ by left multiplication. This semigroup is left Ore and its enveloping group of left quotients is given by $Q(R) \rtimes Q(R)\reg$, where $Q(R)$ is the quotient field of $R$. In the case where $R$ is the ring of integers in a number field, the semigroup C*-algebra of $R \rtimes R\reg$ has been studied intensively in \cite{C-D-L}.

Let us now assume that $R \rtimes R\reg$ satisfies the condition that $\cJ$ is independent. We then observe that since $Q(R) \rtimes Q(R)\reg$ is solvable, the semigroup C*-algebra $C^*(R \rtimes R\reg)$ is nuclear, and full and reduced versions coincide. The boundary quotient of $C^*(R \rtimes R\reg)$ is canonically isomorphic to the ring C*-algebra $\fA[R]$ introduced in \cite{Li1} (compare also \cite{Sun1}, \cite{Sun2} for concrete examples). Moreover, in this case, the group $G_0$ from Proposition~\ref{top-free} coincides with the group of invertible elements in $R \rtimes R\reg$, i.e. $G_0 = R \rtimes R^*$ where $R^*$ is the group of units of $R$. If $R$ is not a field, then it is easy to see that $R \rtimes R^*$ acts topologically freely on $(\partial \Omega) \cdot G$. And by our assumption that $R \neq \gekl{0}$, $R \rtimes R\reg$ is not trivial. Therefore, we can again apply Corollary~\ref{pisun} and deduce that the boundary quotient of $C^*(R \rtimes R\reg)$ is a UCT Kirchberg algebra. As this boundary quotient is nothing else but $\fA[R]$, we have reproven \cite{Li1}, Corollary~8 (for $\cF = \emptyset$).

\section{Open questions and future research}
\label{future}

Of course, one obvious question is how restrictive our assumptions are. It would be interesting to see which semigroups have independent constructible right ideals, and when the Toeplitz condition is satisfied. Is there an intrinsic characterization in terms of the semigroup when a semigroup embeds into a group such that the Toeplitz condition holds? In this context, it would certainly be desirable to study more examples.

In this paper, we have only considered the case of subsemigroups of groups, and one might wonder what to do in the general case of left cancellative semigroups. Recent work in \cite{Nor} and also our results in \S~\ref{invsemi} and \S~\ref{cropro-partialauto} suggest that one should look at left inverse hulls.

One could also try to interprete our results in terms of geometric group theory: Given a subsemigroup $P$ of a group $G$, what is the relationship between nuclearity of the semigroup C*-algebra(s) of $P$ and exactness of $G$? Of course, it would be necessary to impose conditions on $\ping$. Otherwise, one could take the trivial subsemigroup, and the corresponding semigroup C*-algebra is always nuclear. This just reflects the fact that every group acts amenably on itself. But if one asks for the condition that $P$ generates $G$, the problem of relating nuclearity of $C^*_s(P)$ and exactness of $G$ maybe becomes more interesting.

Our main result on nuclearity of semigroup C*-algebras tells us that nuclearity implies faithfulness of the left regular representation. A natural question would be: What about the converse?

One could also study semigroup C*-algebras and their ideals and quotients from the perspective of classification. An interesting question in this context would be which UCT Kirchberg algebras arise as the boundary quotients of semigroup C*-algebras.

\end{document}